\documentclass[twocolumn,twoside]{IEEEtran}

\usepackage{amsmath,epsf}
\numberwithin{equation}{section} % in amsmath

 \newtheorem{lemma}{Lemma}[section]
 \newtheorem{theorem}[lemma]{Theorem}
 \newtheorem{proposition}[lemma]{Proposition}

 \newtheorem{corollary}[lemma]{Corollary}

 \newtheorem{definition}[lemma]{Definition}

 \newtheorem{notation}[lemma]{Notation}

 \newtheorem{ex}[lemma]{Example}
\newenvironment{example}{\begin{ex}}{\hspace*{\fill}$\diamondsuit$\end{ex}}

 \newtheorem{rem}[lemma]{Remark}
\newenvironment{remark}{\begin{rem}}{\hspace*{\fill}$\diamondsuit$\end{rem}}
\newcommand{\lea}{\stackrel{{}_+}{<}}
\newcommand{\gea}{\stackrel{{}_+}{>}}
\newcommand{\eqa}{\stackrel{{}_+}{=}}

\newcommand{\lem}{\stackrel{{}_*}{<}}
\newcommand{\gem}{\stackrel{{}_*}{>}}
\newcommand{\eqm}{\stackrel{{}_*}{=}}

\newcommand{\txt}[1]{\text{\rmfamily\mdseries\upshape{#1}}}
\newcommand{\impl}{\txt{impl}}
\newcommand{\expl}{\txt{expl}}

 \newcommand\setof[1]{\mathopen\{\,#1\,\mathclose\}}
 \newcommand\cei[1]{{\lceil #1\rceil}}

 \newcommand \ag{\alpha}
 \newcommand \bg{\beta}

 \renewcommand{\d}{\delta}
 \newcommand{\m}{{\mathbf{m}}}
 \newcommand{\h}{{h}}

\begin{document}

\title{Algorithmic Statistics}
\author{P\'eter G\'acs, 
%\it Senior Member, IEEE\rm, 
John T. Tromp, and
Paul M.B. Vit\'anyi\thanks{Manuscript 
received June, 2000; revised April 2001.
Part of this work was done during P. G\'acs' stay at CWI.
His work was supported in part by NSF  and by NWO under Grant 
B 62-551.
The work of J.T. Tromp was supported in part
 by NWO under Grant 612.015.001 and by the
EU fifth framework project QAIP, IST--1999--11234, the NoE QUIPROCONE IST--1999--29064,
the ESF QiT Programmme, and the EU Fourth Framework BRA
 NeuroCOLT II Working Group
EP 27150.
The work of P.M.B. Vit\'anyi was supported in part by the 
EU fifth framework project QAIP, IST--1999--11234, the NoE QUIPROCONE IST--1999--29064,
the ESF QiT Programmme, and the EU Fourth Framework BRA
 NeuroCOLT II Working Group
EP 27150.
A preliminary version of part of this work
was published in {\em Proc. 11th Algorithmic Learning Theory Conf.},
{\protect \cite{GTV00}}, and was
archived as http://xxx.lanl.gov/abs/math.PR/0006233.
\newline
P. G\'acs 
is with the Computer Science Department, Boston University, Boston MA 02215,
U.S.A. Email: {\tt gacs@bu.edu}.  
\newline
J.T. Tromp is with CWI, Kruislaan 413, 
1098 SJ Amsterdam, The Netherlands. Email:
{\tt John.Tromp@cwi.nl}.
\newline
P.M.B. Vit\'anyi is with CWI, Kruislaan 413, 
1098 SJ Amsterdam, The Netherlands. Email:
He also has an appointment at the University of Amsterdam.
Email: {\tt Paul.Vitanyi@cwi.nl}.}
}
\markboth{IEEE Transactions on Information Theory, VOL. XX, NO Y, MONTH 2001}{G\'acs, Tromp, and Vit\'anyi: Algorithmic Statistics}

\maketitle
 
\begin{abstract}
  While Kolmogorov complexity is the accepted absolute measure of information
content of an individual finite object, a
similarly absolute notion is needed for the relation between
an individual data sample and an individual model summarizing
the information in the data, for example, 
a finite set (or probability distribution) 
where the data sample typically came from. 
The statistical theory based on such relations between individual objects
can be called algorithmic statistics,
in contrast
to classical statistical theory that deals with relations between probabilistic
ensembles.
%Since the algorithmic theory deals with individual objects and not
%with averages over ensembles of objects it is surprising
%that similar properties hold, albeit sometimes in weaker form.
%We first recall the notion of algorithmic mutual
%information between individual objects and show
%that this information cannot be increased by algorithmic 
%or probabilistic means (as is the case with probabilistic mutual information).
We develop the algorithmic theory of statistic, sufficient
statistic, and minimal sufficient statistic. This theory 
is based on two-part codes
consisting of the code for the statistic (the model summarizing the
regularity, the meaningful information, in the data) 
and the model-to-data code.
In contrast to the situation in
probabilistic statistical theory, the algorithmic
relation of (minimal) sufficiency is an absolute
relation between the individual model and the individual data sample.
We distinguish implicit and explicit descriptions of the models. 
We give characterizations
of algorithmic (Kolmogorov) 
minimal sufficient statistic for all data samples for both
description modes---in the explicit mode under some constraints.
We also strengthen and elaborate earlier results on
the ``Kolmogorov structure function'' and 
``absolutely non-stochastic objects''---those rare objects
for which the simplest models that summarize their
relevant information (minimal sufficient statistics) are at
least as complex as the objects themselves.
We demonstrate a close relation between the probabilistic notions and
the algorithmic ones: (i) in both cases there is an ``information non-increase''
law; (ii) it is shown that a function is a probabilistic sufficient
statistic iff it is with high probability
(in an appropriate sense) an algorithmic sufficient statistic. 
\end{abstract}

\begin{keywords}
Algorithmic information theory;
description format, explicit, implicit;
foundations of statistics;
Kolmogorov complexity;
minimal sufficient statistic, algorithmic;
mutual information, algorithmic;
nonstochastic objects;
sufficient statistic, algorithmic;
two-part codes.
\end{keywords}

\section{Introduction}
\PARstart{S}{tatistical} theory
ideally considers the following problem:
Given a data sample and a family of models (hypotheses), select
the model that produced the data.
But {\em a priori} it is possible that
the data is atypical for the model that actually produced it,
or that the true model is not present in the considered model class. 
Therefore we have to relax our requirements.
%Meaningful induction is possible only by ignoring this possibility.
%Strictly speaking, selection
If selection of a ``true'' model cannot be guaranteed by any method, 
then as next best choice
``modeling the data'' as well as possible irrespective 
of truth and falsehood of
the resulting model may be more appropriate.
Thus, we change ``true'' to ``as well as possible.'' 
The latter we take to mean that the model expresses all significant
regularity present in the data. 
The general setting is as follows: We carry
out a probabilistic experiment of which the outcomes are governed
by an unknown probability distribution $P$.
Suppose we obtain as outcome the
data sample $x$. Given 
$x$, we want to recover the distribution $P$. For certain reasons
we can choose a distribution from a set of acceptable distributions
only (which may or may not contain $P$). Intuitively, our selection criteria
are that (i) $x$ should be a ``typical'' outcome of the distribution
selected, and (ii) the selected distribution has a ``simple'' description.
We need to make the meaning of ``typical'' and ``simple'' rigorous
and balance the requirements (i) and (ii). In probabilistic statistics
one analyzes the average-case performance of the selection process.
For traditional problems, dealing with frequencies over small
sample spaces, this approach is appropriate. But for 
current novel applications, average relations are often 
irrelevant, since the part of the support of the probability
density function that will ever be observed has about zero
measure. This is the case in, for example, complex video and sound
analysis.
There arises the problem that for individual cases the selection 
performance may be bad although the performance is good on average.
We embark on a systematic study of model selection where
the performance is related to the individual data sample and the individual 
model selected. It turns out to be more straightforward to investigate
models that are finite sets first, and then generalize the results
to models that are probability distributions. To simplify matters,
and because all discrete data
can be binary coded, we consider only data samples
that are finite binary strings.

This paper is one of a triad of papers dealing with the
best individual model for individual data: The present paper supplies the
basic theoretical underpinning by way of two-part codes, \cite{ViLi99}
derives ideal versions of applied methods (MDL) inspired by the theory, and
\cite{GLV00} treats experimental applications thereof.

{\bf Probabilistic Statistics:}
In ordinary statistical theory one proceeds as follows, see for example
\cite{CT91}:
Suppose two discrete random variables $X,Y$ have a joint probability
mass function $p(x,y)$ and marginal probability mass functions 
$p_1(x) = \sum_y p(x,y)$
and $p_2(y) = \sum_x p(x,y)$. Then the (probabilistic)
 {\em mutual information} $I(X;Y)$ 
between the joint distribution and the product distribution $p_1(x)p_2(y)$ is
defined by:
\begin{equation}\label{eq.mutinfprob}
 I(X;Y) = \sum_x \sum_y p(x,y) \log \frac{p(x,y)}{p_1(x)p_2(y)} ,
\end{equation}
where ``$\log$'' denotes the binary logarithm.
Consider a probabilistic ensemble of models,
say a family of probability mass functions $\{f_{\theta} \}$
indexed by $\theta$, together with a distribution $p_1$ over $\theta$.
This way we have a random variable $\Theta$ with outcomes in $\{f_{\theta} \}$
and a random variable $D$ with outcomes in the union of domains of $f_{\theta}$,
and $p(\theta,d) = p_1 (\theta) f_{\theta}(d)$. 
Every function $T(D)$ of a data sample $D$---like the sample
mean or the sample variance---is called
a {\em statistic} of $D$. 
A statistic $T(D)$ is called {\em sufficient} if the probabilistic mutual 
information 
\begin{equation}\label{eq.suffstatprob}
I(\Theta ; D) = I( \Theta ; T(D))
\end{equation}
 for all distributions
of $\theta$.
Hence, the mutual information between  parameter and data sample
random variables is invariant under taking sufficient statistic and vice versa.
That is to say, a statistic $T(D)$ is called sufficient
for $\Theta$ if it contains all the
information in $D$ about $\Theta$.
For example, consider $n$ tosses of a coin with unknown bias $\theta$
with outcome $D=d_1 d_2 \ldots d_n$ where $d_i \in \{0,1\}$ ($1 \leq i \leq n$).
Given $n$, the number of outcomes ``1'' is a sufficient statistic for
$\Theta$: the statistic $T(D) = s = \sum_{i=1}^n d_i$. Given $T$,
all sequences with $s$ ``1''s are equally likely independent
of parameter $\theta$: Given $s$, if $d$ is an outcome of $n$ coin tosses
and $T(D)=s$ then $\Pr(d \mid T(D)=s) = {n \choose s}^{-1}$
and  $\Pr(d \mid T(D) \neq s) = 0$. This can be
shown to imply (\ref{eq.suffstatprob})
and therefore $T$ is a sufficient statistic for $\Theta$.
According to Fisher~\cite{Fi22}:
       ``The statistic chosen should summarise the whole of the relevant
information supplied by the sample. This may be called
       the Criterion of Sufficiency $\ldots$ 
In the case of the normal curve
of distribution it is evident that the second moment is a
       sufficient statistic for estimating the standard deviation.''
Note that one cannot improve on sufficiency:
 for every (possibly randomized) function $T$
we have 
\begin{equation}\label{eq.infnonincrprob}
I(\Theta ; D) \geq I( \Theta ; T(D)),
\end{equation}
that is, mutual information
cannot be increased by processing the data sample in any way.

A sufficient statistic may contain information
that is not relevant: for a normal distribution the sample mean
is a sufficient statistic, but the pair of functions
which give the 
mean of the even-numbered samples and the odd-numbered samples
respectively, is also a sufficient statistic.
A statistic $T(D)$ is a {\em minimal} sufficient statistic
with respect to an indexed
model family $\{f_{\theta}\}$, if it is a
function of all other sufficient statistics: it contains no
irrelevant information and maximally compresses the information about
the model ensemble.
As it happens, for the family of normal distributions
the sample mean is a minimal sufficient statistic, but the 
sufficient statistic consisting of the mean of the even samples
in combination with the mean of the odd samples is not minimal.
All these notions and laws are probabilistic: they hold
in an average sense.

{\bf Kolmogorov Complexity:}
  We write {\em string} to mean a finite binary sequence.
  Other finite objects can be encoded into strings in natural
ways.
  The Kolmogorov complexity, or algorithmic entropy, $K(x)$ of a
string $x$ is the length of a shortest binary program to compute
$x$ on a universal computer (such as a universal Turing machine).
  Intuitively, $K(x)$ represents the minimal amount of information
required to generate $x$ by any effective process, \cite{Ko65}.
  The conditional Kolmogorov complexity $K(x \mid y)$ of $x$ relative to
$y$ is defined similarly as the length of a shortest program
to compute $x$ if $y$ is furnished as an auxiliary input to the
computation. This conditional definition
requires a warning since different authors use the same notation
but mean different things. In \cite{Ch75} the author writes
\label{fo.gacs}
``$K(x \mid y)$'' to actually mean ``$K(x \mid y, K(y))$,''
notationally hiding the intended supplementary 
auxiliary information ``$K(y)$.''
This abuse of notation has the additional handicap
that no obvious notation is left to express ``$K(x \mid y)$''
meaning that just ``$y$'' is given in the conditional.
As it happens,
``$y, K(y)$'' represents more information than just ``$y$''. For example,
$K(K(y) \mid y)$ can be almost as large as $\log K(y)$ by a result
in \cite{Ga74}: 
For $l(y)=n$ it has an upper bound of $\log n$ for all $y$, and
for some $y$'s it has a lower bound of $\log n - \log \log n$.
In fact, this result quantifies the undecidability
of the halting problem for Turing machines---for example,
if $K(K(y) \mid y) = O(1)$ for all $y$, then the halting problem
can be shown to be decidable. This is known to be false.
It is customary, \cite{Le74,Ga74,LiVi97}, to write explicitly
``$K(x \mid y)$'' and ``$K(x \mid y, K(y))$''.
Even though the difference between these two quantities is not very
large,
these small
differences do matter in the sequel. In fact, not only the precise
information itself in the conditional, but also the way it is 
represented, is crucial, see 
Subsection~\ref{ss.exim}.

  The functions $K( \cdot)$ and $K( \cdot \mid  \cdot)$,
though defined in terms of a
particular machine model, are machine-independent up to an additive
constant
 and acquire an asymptotically universal and absolute character
through Church's thesis, from the ability of universal machines to
simulate one another and execute any effective process.
  The Kolmogorov complexity of a string can be viewed as an absolute
and objective quantification of the amount of information in it.
   This leads to a theory of {\em absolute} information {\em contents}
of {\em individual} objects in contrast to classical information theory
which deals with {\em average} information {\em to communicate}
objects produced by a {\em random source}.
  Since the former theory is much more precise, it is surprising that
analogs of theorems in classical information theory hold for
Kolmogorov complexity, be it in somewhat weaker form. Here our aim is
to provide a similarly absolute notion for individual ``sufficient statistic''
and related notions borrowed from probabilistic statistics. 

{\bf Two-part codes:}
The prefix-code of the shortest effective descriptions
gives an expected code word length close to the entropy
and also compresses the regular objects until all regularity is
squeezed out. All shortest effective descriptions are
completely random themselves, without any regularity whatsoever.
The idea of a two-part code for a body of data $d$
is natural from the perspective of Kolmogorov
complexity.
If $d$ does not contain any regularity at all, then it consists
of purely random data and the model is precisely that.
Assume that the body of data $d$ contains regularity.
With help of a description of the regularity (a model) we can
describe the data compactly. Assuming that the regularity can be represented
in an effective manner (that is, by a Turing machine),
we encode the data as a program for that machine. Squeezing
all effective regularity out of the data, we end up
with a Turing machine representing the meaningful regular
information in the data together with a program for
that Turing machine representing the remaining
meaningless randomness of the data. 
However, in general there are many ways
to make the division into meaningful information and remaining
random information. In a painting
the represented image, the brush strokes, or even
finer detail can be the relevant information,
depending on what we are interested in. What we require is
a rigorous mathematical condition to force a sensible division
of the information at hand into a meaningful part
and a meaningless part.

{\bf Algorithmic Statistics:}
The two-part code approach leads to a more general algorithmic
approach to statistics.
The algorithmic statistician's task is to
 select a model
(described possibly by a probability distribution)
for which the data is typical.
In a two-part description, 
 we describe such a model and then identify
the data within the set of the typical outcomes.
The best models make the two-part description
as concise as the best
one-part description of the data.
A description of
such a model is an algorithmic sufficient statistic
since it summarizes all relevant properties of the data.
Among the algorithmic sufficient statistics,
the simplest one (an algorithmic minimal sufficient statistic) 
is best in accordance
with Ockham's Razor since it summarizes the relevant
properties of the data as concisely as possible.
In probabilistic data or data 
subject to noise this involves separating regularity (structure)
in the data from random effects.

In a restricted setting where the models
are finite sets a way to proceed 
was suggested by Kolmogorov, attribution in \cite{ShenStoch83,Co85,CT91}. 
Given data $d$,
the goal is to identify
the ``most likely'' finite set $S$ of which $d$ is a ``typical'' element.
Finding a set of which the data is typical is 
reminiscent of selecting the appropriate magnification of a microscope to
bring the studied specimen optimally in focus.
For this purpose we consider sets $S$ such that
$d \in S$ and we represent $S$ by the {\em shortest} program $S^*$ that
computes the characteristic function of $S$. 
The shortest program $S^*$ that computes 
a finite set $S$ containing $d$,
such that
the two-part description consisting of $S^*$ and $\log |S|$ is as
as short as the shortest {\em single} program that computes $d$ without
input, is called
an {\em algorithmic 
sufficient statistic}\footnote{It is also called the Kolmogorov 
sufficient statistic.}
This definition is non-vacuous since
there does exist a two-part code (based on the model $S_d = \{d\}$)
that is as concise as the shortest single code.
The description of $d$ given $S^*$
cannot be significantly shorter than $\log |S|$. By
the theory of Martin-L\"of randomness \cite{ML66} this means
that $d$ is a ``typical'' element of $S$.
In general there can be many algorithmic sufficient statistics
for data $d$; a shortest among them is called an {\em algorithmic minimal
sufficient statistic}.
Note that there can be possibly more than one
algorithmic minimal sufficient statistic; they are
defined by, but not generally computable from, the data.

In probabilistic statistics the notion of sufficient statistic
(\ref{eq.suffstatprob}) is an average notion invariant under all probability
distributions over the family of indexed models. If a statistic
is not thus invariant, it is not sufficient.
In contrast, in the algorithmic case we investigate 
the relation between
the data and an individual model and therefore
a probability distribution
over the models is irrelevant.
It is technically convenient
to initially consider the simple model class of finite sets to 
obtain our results. It then turns out that it is relatively easy to
generalize everything to the model class of computable
probability distributions. That class is very large indeed: perhaps
it contains every distribution that has  ever been considered in
statistics and probability theory, as long as the parameters
are computable numbers---for example rational numbers. Thus the
results are of great generality; indeed, they are so general
that further development of the theory must be aimed at restrictions
on this model class, see the discussion about applicability
in Section~\ref{s.discussion}.  
The theory concerning
the statistics of individual data samples and models
one may call {\em algorithmic statistics}.

{\bf Background and Related Work:}
At a Tallinn conference in 1973,
A.N. Kolmogorov formulated 
the approach to an individual data to model relation,
based on a two-part code separating
the {\em structure} of a string from meaningless {\em random} features,
rigorously in terms of 
Kolmogorov complexity (attribution by  \cite{ShenStoch83,Co85}). 
Cover~\cite{Co85,CT91} interpreted this approach as
a (sufficient) statistic. The ``statistic'' of the data 
is expressed as a finite set of which the data is a ``typical''
member.
Following Shen~\cite{ShenStoch83} 
(see also~\cite{Vy87,ShenStoch99,ViLi99}), this can be generalized
to computable probability mass functions for which the data
is ``typical.''
Related aspects of ``randomness deficiency'' (formally defined later
in (\ref{eq.randdef}))
were formulated in~\cite{Ko83,KU88} and
studied in~\cite{ShenStoch83,Vy87}.
Algorithmic mutual information, and the associated non-increase law,
were studied in \cite{Le74,LevinRandCons84}.
Despite its evident epistemological
prominence in the theory of hypothesis selection
and prediction, only selected
aspects of the algorithmic sufficient statistic
have been studied before, for example as related to the
``Kolmogorov structure function''~\cite{ShenStoch83,Co85},
 and ``absolutely non-stochastic objects''~\cite{ShenStoch83,Vy87,ShenStoch99,Vy99}, notions also defined or suggested by Kolmogorov
at the mentioned meeting. 
This work primarily studies quantification of the ``non-sufficiency''
of an algorithmic statistic, when the latter is restricted in
complexity, rather than necessary and sufficient
conditions for the existence of an algorithmic sufficient statistic itself.
These references obtain results for plain 
Kolmogorov complexity (sometimes length-conditional)
up to a  logarithmic error term. 
Especially for regular data that have low Kolmogorov
complexity with respect to their length, this logarithmic error term
may dominate the remaining terms and eliminate
all significance.
Since it is precisely the regular data that
one wants to assess the  meaning of, a more precise analysis as we 
provide is required.
Here we use prefix complexity
%pg: (not length-conditional) 
% In my part, it is length-conditional.
to unravel the nature of a sufficient statistic.
The excellent papers of Shen \cite{ShenStoch83,ShenStoch99}
contain the major previous results related to this work (although
\cite{ShenStoch99} is independent). While previous work and the present
paper consider an algorithmic statistic that is either a finite set
or a computable probability mass function, the most general algorithmic
statistic is a recursive function. In \cite{AFV} the present work
is generalized accordingly, see the summary in Section~\ref{s.discussion}.   

For the relation with
inductive reasoning according to minimum description length principle see
~\cite{ViLi99}.
The entire approach is based on
Kolmogorov complexity (also known as algorithmic information
theory). Historically, the idea of assigning to each object a
probability
consisting of the summed negative exponentials of the lengths
of all programs
computing the object, was first proposed
by Solomonoff \cite{So64}. Then, the shorter programs
contribute more probability than the longer ones.
His aim, ultimately successful in terms of theory (see \cite{LiVi97})
 and as inspiration
for developing applied versions \cite{BRY},  was to develop
a general prediction method.
Kolmogorov  \cite{Ko65} introduced the complexity proper. 
The prefix-version of Kolmogorov complexity used in this paper
was introduced in \cite{Le74} and also
treated later in \cite{Ch75}.
For a textbook on Kolmogorov complexity,
its mathematical theory, and its application to induction, see~\cite{LiVi97}.
We give a definition (attributed to Kolmogorov)
and results from \cite{ShenStoch83} that are useful later:

\begin{definition}
Let $\alpha$ and $\beta$ be natural numbers.  A finite binary string
$x$ is called {\em $(\alpha, \beta )$-stochastic} if there exists a finite
set $S \subseteq \{0,1\}^*$ such that
\begin{equation}\label{def.stoch}
 x \in S, \; \; K(S) \leq \alpha , \; \; K(x) \geq \log |S| - \beta ;
\end{equation}
where $|S|$ denotes the cardinality of $S$, and $K(\cdot)$
the (prefix-) Kolmogorov complexity. As usual, ``$\log$'' denotes the binary
logarithm.
\end{definition}

The first inequality with small $\alpha$ means that $S$ is ``simple'';
the second inequality with $\beta$ is small means that $x$ is
``in general position'' in $S$. Indeed, if $x$ had any
special property $p$ that was shared by only a small subset $Q$
of $S$, 
then this property could be used to single out and enumerate those
elements and subsequently indicate $x$ by its index in the enumeration.
Altogether, this would show $K(x) \leq K(p) + \log |Q|$, which,
for simple $p$ and small $Q$ would be much lower than $\log |S|$.
A similar notion for computable probability distributions
is as follows: 
Let $\alpha$ and $\beta$ be natural numbers.  A finite binary string
$x$ is  called {\em $(\alpha, \beta)$-quasistochastic} if there exists
a computable probability distribution $P$ such that
\begin{equation}\label{def.quasistoch}
P(x) > 0, \; \; K(P) \leq \alpha, \; \; K(x) \geq - \log P(x) - \beta .
\end{equation}
\begin{proposition}\label{prop.1}
There exist constants $c$ and $C$,
 such that for every natural number $n$
and every finite binary string $x$ of length $n$:

(a) if $x$ is $(\alpha, \beta)$-stochastic, then $x$ is 
$(\alpha + c, \beta)$-quasistochastic; and

(b) if $x$ is $(\alpha, \beta)$-quasistochastic and the length of $x$
is less than $n$, then $x$ is $(\alpha  + c \log n, \beta + C)$-stochastic.
\end{proposition}

\begin{proposition}\label{prop.2}
(a)  There exists a constant $C$ such that, 
for every natural number $n$ and every $\alpha$ and $\beta$
with $\alpha \geq \log n + C$ and $\alpha + \beta \geq n + 4 \log n + C$,
all strings of length less than $n$ are $(\alpha, \beta)$-stochastic.

(b) 
There exists a constant $C$ such that, 
for every natural number $n$ and every $\alpha$ and $\beta$
with $2 \alpha + \beta < n - 6 \log n - C$, 
there exist strings $x$ of length less than
$n$ that are not $(\alpha, \beta)$-stochastic.
\end{proposition}
Note that if we take $\alpha = \beta$ then, for some
boundary in between $\frac{1}{3}n$ and $\frac{1}{2}n$, 
the last non-$(\alpha, \beta)$-stochastic elements disappear if
the complexity constraints are sufficiently
relaxed by having $\alpha, \beta$ exceed this boundary.

{\bf Outline of this Work:}
First, we obtain a new Kolmogorov complexity ``triangle''
inequality that is useful in the later parts of the paper.
We define algorithmic mutual information between two
individual objects (in contrast to the probabilistic notion of
mutual information that deals with random variables).
We show that for every computable distribution
associated with the random variables,
the expectation of the algorithmic mutual information
equals the probabilistic mutual information up to an additive
constant that depends on the complexity of the distribution.
It is known that in the probabilistic setting the mutual information
(an average notion)
cannot be increased by algorithmic processing. We give a new proof that this
also holds in the individual setting. 

We define notions of ``typicality'' and ``optimality'' of sets in
relation to the given data $x$. Denote the shortest program for a
finite set $S$ by $S^*$ (if there is more than one shortest program
$S^*$ is the first one in the standard effective enumeration).
``Typicality'' is a reciprocal relation:
A set $S$ is ``typical'' with respect to $x$  if 
$x$ is an element of $S$ that is ``typical''
in the sense of having small {\em randomness
deficiency} $\delta_S^* (x) = \log |S|-K(x|S^*)$ (see
definition (\ref{eq.randdef}) and discussion).
That is, $x$ has about maximal Kolmogorov complexity
in the set, because it can always be identified by its position
in an enumeration of $S$ in $\log |S|$ bits. 
Every description of a ``typical'' set for the data is an 
algorithmic statistic.

A set $S$ is ``optimal''  if the best two-part description consisting
of a description of $S$ and a 
%pg:
straightforward %
description of $x$ as an element of $S$
%pg: 
by an index of size $\log |S|$
%pg: Without this, it was not clear to me.
is as concise as the shortest one-part description of $x$.
This implies 
that optimal sets are typical sets.
Descriptions of such optimal sets
are algorithmic sufficient statistics, and a shortest description
among them is an algorithmic minimal sufficient statistic.
The mode of
description plays a major role in this.
We distinguish between
``explicit'' descriptions and ``implicit'' descriptions---that are
introduced in this paper as a proper restriction 
on the recursive enumeration based description mode.
We establish range constraints of cardinality and complexity imposed
by implicit (and hence explicit) descriptions for typical and optimal
sets, and exhibit a concrete algorithmic minimal
sufficient statistic for implicit description mode.
It turns out that
only the complexity of the data sample $x$ is relevant for this
implicit algorithmic minimal sufficient statistic.
Subsequently we 
exhibit explicit algorithmic sufficient
statistics, and an explicit minimal algorithmic (near-)sufficient statistic.
For explicit descriptions it turns out that certain other
aspects of $x$  (its enumeration rank) apart from its complexity
are a major determinant for the cardinality and complexity 
of that statistic. 
It is convenient at this point
to introduce some notation:

\begin{notation}
\rm
  From now on, we will denote by $\lea$ an inequality to within an
additive constant, and by $\eqa$ the situation when both $\lea$ and
$\gea$ hold.
  We will also use $\lem$ to denote an inequality to within an
multiplicative constant factor, and $\eqm$ to denote
the situation when both $\lem$ and
$\gem$ hold.
\end{notation}

Let us contrast our approach with the 
one in \cite{ShenStoch83}. The comparable case there,
by (\ref{def.stoch}), is that $x$ is 
$(\alpha, \beta)$-stochastic with $\beta = 0$ and $\alpha$ minimal. 
Then, $K(x) \geq \log |S|$
for a set $S$ of Kolmogorov complexity $ \alpha$.
But, if $S$ is  optimal for $x$, then, as we formally define it later
(\ref{eq.optim}),
$K(x) \eqa K(S)+ \log |S|$. 
That is (\ref{def.stoch}) holds with $\beta \eqa - K(S)$. In contrast,
for $\beta = 0$ we must have $K(S) \eqa 0$ for typicality.
In short, optimality of $S$ with repect to $x$
corresponds to (\ref{def.stoch})
by dropping the second item and replacing the third item 
by $K(x) \eqa \log |S| + K(S)$. ``Minimality'' of the algorithmic
sufficient statistic $S^*$ (the shortest program for $S$)
corresponds to choosing $S$ with minimal $K(S)$ in this equation. 
This is equivalent to (\ref{def.stoch}) with inequalities replaced
by equalities and $K(S)= \alpha = - \beta$.

We consider the functions related to
$(\alpha, \beta)$-stochasticity, and improve Shen's
result on maximally non-stochastic objects.
In particular, we show that for every $n$ there are
objects $x$ of length $n$  with complexity $K(x \mid n)$ about $n$  
such that every explicit algorithmic  sufficient statistic
for $x$  has complexity about $n$ ($\{x\}$ is such a statistic).
This is the best possible.
In Section~\ref{s.prob}, we
generalize the entire treatment to probability density distributions.
In Section~\ref{sect.formanal} we connect the algorithmic and
probabilistic approaches:
While previous authors have used the name ``Kolmogorov sufficient statistic''
because the model appears to summarize the relevant information in the data
in analogy of what the classic sufficient statistic 
does in a probabilistic sense, a formal justification has been lacking.
We give the formal relation between the
algorithmic approach to sufficient statistic and the probabilistic
approach: A function is a probabilistic sufficient statistic iff 
it is with high probability an algorithmic
$\theta$-sufficient statistic, where an algorithmic
sufficient statistic is {\em $\theta$-sufficient} if it 
satisfies also the sufficiency criterion conditionalized on $\theta$.

 \section{Kolmogorov Complexity} 
We give some definitions to establish notation.
For introduction, details, and proofs, see \cite{LiVi97}.
We write {\em string} to mean a finite binary string.
  Other finite objects can be encoded into strings in natural
ways.
  The set of strings is denoted by $\{0,1\}^*$. The {\em length}
of a string $x$ is denoted by $l(x)$, distinguishing it
from the {\em cardinality} $|S|$ of a finite set $S$.

Let $x,y,z \in \mathcal{N}$, where
$\mathcal{N}$ denotes the natural
numbers.
Identify
$\mathcal{N}$ and $\{0,1\}^*$ according to the
correspondence
 \[
 (0, \epsilon ), (1,0), (2,1), (3,00), (4,01), \ldots . 
 \]
Here $\epsilon$ denotes the {\em empty word} `' with no letters.
The {\em length} $l(x)$ of $x$ is the number of bits
in the binary string $x$. For example,
$l(010)=3$ and $l(\epsilon)=0$.

The emphasis is on binary sequences only for convenience;
observations in any alphabet can be so encoded in a way
that is `theory neutral'.

A binary string $x$
is a {\em proper prefix} of a binary string $y$
if we can write $y=xz$ for $z \neq \epsilon$.
 A set $\{x,y, \ldots \} \subseteq \{0,1\}^*$
is {\em prefix-free} if for any pair of distinct
elements in the set neither is a proper prefix of the other.
A prefix-free set is also called a {\em prefix code}.
Each binary string $x=x_1 x_2 \ldots x_n$ has a
special type of prefix code, called a
{\em self-delimiting code},
\[ \bar x = 1^n 0 x_1x_2 \ldots x_n  .\]
This code is self-delimiting because we can determine where the
code word $\bar x$ ends by reading it from left to right without
backing up. Using this code we define
the standard self-delimiting code for $x$ to be
$x'=\overline{l(x)}x$. It is easy to check that
$l(\bar x ) = 2 n+1$ and $l(x')=n+2 \log n+1$.

Let $\langle \cdot ,\cdot \rangle$ be a standard one-one mapping
from $\mathcal{N} \times \mathcal{N}$
to $\mathcal{N}$, for technical reasons chosen such that
$l(\langle x, y \rangle) = l(y) + l(x) + 2 l( l(x)) +1$, for example
$\langle x, y \rangle = x'y =
1^{l(l(x))}0l(x)xy$.
This can be iterated to
$\langle  \langle \cdot , \cdot \rangle , \cdot \rangle$.

  The {\em prefix Kolmogorov complexity}, 
or algorithmic entropy, $K(x)$ of a
string $x$ is the length of a shortest binary program to compute
$x$ on a universal computer (such as a universal Turing machine).
For technical reasons we require that the universal machine has
the property that no halting program is a proper prefix of another
halting program.
  Intuitively, $K(x)$ represents the minimal amount of information
required to generate $x$ by any effective process.
We denote the {\em shortest program} for $x$ by $x^*$; then
$K(x)= l(x^*)$. 
(Actually, $x^*$ is the first shortest program for $x$ in
an appropriate standard enumeration of all programs for $x$
such as the halting order.)
  The conditional Kolmogorov complexity $K(x \mid y)$ of $x$ relative to
$y$ is defined similarly as the length of a shortest program
to compute $x$ if $y$ is furnished as an auxiliary input to the
computation. We often use
$K(x \mid y^*)$, or, equivalently, $K(x \mid y,K(y))$
(trivially $y^*$ contains the same information
as the $y,K(y)$). Note that ``$y$'' in the conditional
is just the information about $y$ and apart from
this does not contain information
about $y^*$ or $K(y)$. For this work the difference is
crucial, see the comment in Section~\ref{fo.gacs}.
  %The functions $K( \cdot)$ and $K( \cdot \mid  \cdot)$,
%though defined in terms of a
%particular universal machine model, are machine-independent up to an 
%additive constant
 %and acquire an asymptotically universal and absolute character
%through Church's thesis, from the ability of universal machines to
%simulate one another and execute any effective process.
  %The Kolmogorov complexity of a string can be viewed as an absolute
%and objective quantification of the amount of information in it.
  %% \footnote{This leads to a theory of {\em absolute} information {\em contents}
%%of {\em individual} objects in contrast to classical information theory
%%which deals with {\em average} information {\em to communicate}
%%objects produced by a {\em random source}.
  %%Since the former theory is much more precise, it is surprising that
%%analogons of theorems in classical information theory hold for
%%Kolmogorov complexity, be it in somewhat weaker form.}

\subsection{Additivity of Complexity}

Recall that by definition $K(x,y) = K(\langle x,y \rangle)$. 
Trivially, the symmetry property holds: $K(x,y) \eqa K(y,x)$.
Later we will use many times the ``Additivity of Complexity'' 
property 
 \begin{equation}\label{eq.soi}
  K(x, y) \eqa K(x) + K(y \mid x^*) \eqa K(y) + K(x \mid y^*).
 \end{equation}
This result due to \cite{Ga74} can be found
as Theorem 3.9.1 in~\cite{LiVi97} and has a difficult proof.
It is perhaps instructive to point out that
the version with just $x$ and $y$ in the conditionals doesn't
hold with $\eqa$, but holds up to additive logarithmic terms 
that cannot be eliminated.
The conditional version needs to be treated carefully.
It is
 \begin{equation}\label{eq.soi-cond}
  K(x, y \mid z) \eqa K(x \mid z) + K(y \mid x, K(x \mid z), z).
 \end{equation}
Note that a naive version
 \[
  K(x, y \mid z) \eqa K(x \mid z) + K(y \mid x^{*}, z)
 \]
is incorrect: taking $z = x$, $y = K(x)$,
the left-hand side equals $K(x^{*} \mid x)$, and the right-hand side
equals $K(x \mid x) + K(K(x) \mid x^{*}, x) \eqa 0$.
First, we derive a (to our knowledge) new ``directed triangle inequality''
that is needed later.
 \begin{theorem}\label{lem.magic}
For all $x,y,z$, 
 \[
  K(x \mid y^*) \lea K(x, z \mid y^{*}) \lea K(z \mid y^*) + K(x \mid z^*).
 \]
 \end{theorem}

\begin{proof}
Using~(\ref{eq.soi}), an evident inequality introducing
an auxiliary object $z$, and twice (~\ref{eq.soi}) again:
 \begin{align*}
  K(x, z \mid y^*) &\eqa 
    K(x,y,z) - K(y) 
\\ & \lea K(z) + K(x \mid z^*) + K(y \mid z^*) - K(y)
\\ &\eqa K(y,z) - K(y) + K(x \mid z^*) 
\\ & \eqa K(x \mid z^*) + K(z \mid y^*).
 \end{align*}

\end{proof}

This theorem has bizarre consequences. These  consequences are not
simple unexpected artifacts of our definitions, but, to the contrary,
they show the power and the genuine contribution to our understanding
represented by the deep and important mathematical relation
(\ref{eq.soi}).
 
Denote $k=K(y)$ and substitute $k=z$ and $K(k)=x$
to find the following counterintuitive corollary: To determine the complexity
of the complexity of an object $y$ it suffices to give both $y$ and
the complexity of $y$. This is counterintuitive since in general
we cannot compute the complexity of an object from the object itself;
if we could this would also solve the 
so-called ``halting problem'', \cite{LiVi97}. This noncomputability
can be quantified in terms of $K(K(y) \mid y )$ which can rise to
almost $K(K(y))$ for some $y$---see the related discussion
on notation for conditional complexity in Section~\ref{fo.gacs}. But in the
seemingly similar, but subtly different, setting below it is possible.

\begin{corollary}
As above, let $k$ denote $K(y)$. Then,
$K(K(k) \mid  y,k) \eqa K(K(k) \mid y^*)  \lea K(K(k) \mid k^*)+K(k \mid y,k) \eqa 0$.
We can iterate this idea.
For example, the next step is that
 given $y$ and $K(y)$ we can determine $K(K(K(y)))$
in $O(1)$ bits, 
that is, $K(K(K(k))) \mid y,k) \eqa 0$. 
\end{corollary}

A direct construction works according to the following idea (where
we ignore some important details):
 From $k^*$ one can compute $\langle k, K(k) \rangle$  since
$k^*$ is by definition  the shortest program for $k$
and also by definition $l(k^*)=K(k)$.
Conversely, from $k,K(k)$ one can compute $k^*$: by 
running of all programs of length at most $K(k)$
in dovetailed fashion until the first programme of length $K(k)$
halts with output $k$;
this is $k^*$.
The shortest program that computes the pair $\langle y,k \rangle$
has length $\eqa k$:
We have $K(y,k) \eqa k$ (since the shortest program
$y^*$ for $y$ carries both the information about $y$ and about $k=l(y^*)$).
By~(\ref{eq.soi}) therefore
$K(k)+K(y \mid k, K(k)) \eqa k$. In view of the information
equivalence of $\langle k,K(k) \rangle$ and $k^*$, 
therefore $K(k) + K(y \mid k^*) \eqa k$. Let $r$ be a program
of length $l(r)=K(y \mid k^*)$ that computes $y$ from $k^*$.
Then, since $l(k^*) = K(k)$, there is a shortest
program $y^* =  qk^* r$ for $y$
where $q$ is a  fixed $O(1)$ bit self-delimiting
program that unpacks and uses $k^*$ and $r$ 
to compute $y$.
We are now in the position to show $K(K(k) \mid  y,k) \eqa 0$.
There is a fixed $O(1)$-bit program, that includes knowledge of $q$,
and that enumerates two lists
in parallel, each in dovetailed fashion:
Using $k$ it enumerates a list of 
all programs that compute $k$, including $k^*$.
Given $y$ and $k$ it
enumerates another list of all programs of 
length $k \eqa l(y^*)$ that compute $y$.
One of these programs is $y^* = q k^* r$
that starts with $q k^*$. Since $q$ is known,
this self-delimiting program
 $k^*$, and hence its length $K(k)$, can be found by matching
every element in the $k$-list with the prefixes of every element in 
the $y$ list in enumeration order.

\subsection{Information Non-Increase}

If we want to find an appropriate model fitting the data, then
we are concerned with the information in the data about such models.
Intuitively one feels that the information
in the data about the appropriate model cannot be increased
 by any algorithmic or
probabilistic process.
Here, we rigorously show that this is the case in the algorithmic statistics
setting:
the information in one object about another
cannot be increased by any deterministic algorithmic method
by more than a constant. With added randomization this holds
with overwhelming probability.
We use the triangle inequality of Theorem~\ref{lem.magic} to recall,
and to give possibly new proofs, of this information non-increase;
for more elaborate but hard-to-follow versions see
\cite{Le74,LevinRandCons84}.

We need the following technical concepts.
Let us call a nonnegative 
real function $f(x)$ defined on strings a {\em semimeasure} if 
$\sum_{x} f(x) \le 1$, and a {\em measure} (a probability distribution)
if the sum is 1.
A function $f(x)$ is called {\em lower semicomputable} if there is a 
rational valued computable function $g(n,x)$ such that
$g(n+1,x) \geq g(n,x)$ and $\lim_{n \rightarrow \infty} g(n,x) = f(x)$.
For an {\em upper semicomputable} function $f$ we require 
that $-f$ is lower semicomputable.
It is computable when it is both lower and upper semicomputable.
(A lower semicomputable measure is also computable.)

To define the algorithmic mutual information between
two individual objects $x$ and $y$ with no
probabilities involved, it is instructive to first recall
the probabilistic notion (\ref{eq.mutinfprob})
Rewriting (\ref{eq.mutinfprob})
as 
\[ \sum_x \sum_y p(x,y) [ - \log p(x) - \log p(y) + \log p(x,y) ] , \]
and noting that $- \log p ( s )$ is 
very close to the length of the 
prefix-free Shannon-Fano code for $s$, we are led to the following
definition.
\footnote{The Shannon-Fano code has nearly optimal expected
code length equal to the entropy with 
respect to the distribution of the source \cite{CT91}. However,
the prefix-free code with code word length $K(s)$ has both
about expected optimal code word length and individual optimal
effective code word length, \cite{LiVi97}.}
The
{\em information in  $y$ about $x$}
 is defined as
 \begin{equation}\label{def.mutinf}
   I(y : x) = K(x) - K(x  \mid  y^*) \eqa K(x) + K(y) - K(x, y),
 \end{equation}
where the second equality is a consequence of~(\ref{eq.soi})
and states that this information is symmetrical,
$I(x:y) \eqa I(y:x)$, and therefore we can talk about
{\em mutual information}.\footnote{The notation of the
algorithmic (individual) 
 notion $I(x:y)$ distinguishes it from the probabilistic
(average) notion 
$I(X; Y)$.  We deviate slightly from~\cite{LiVi97}
where $I(y : x)$ is defined as $K(x) - K(x \mid y)$.}
 \begin{remark}\label{rem.cami}
The conditional mutual information is
 \begin{align*}
   I(x : y \mid z) & = K(x \mid z) - K(x \mid y, K(y \mid z), z)
 \\ & \eqa K(x \mid z) + K(y \mid z) - K(x, y \mid z).
 \end{align*}
 \end{remark}
It is important that the expectation of the algorithmic mutual 
information $I(x:y)$ is close to the probabilistic mutual information
$I(X; Y)$---if
this were not the case then the algorithmic notion would not
be a sharpening of the probabilistic notion to individual objects,
but something else.

\begin{lemma}
Given a computable joint probability mass distribution $p(x,y)$ over $(x,y)$
we have 
\begin{align}\label{eq.eqamipmi}
I(X; Y) - K(p) & \lea  \sum_x \sum_y p(x,y) I(x:y) 
\\& \lea I(X;Y) + 2 K(p) , 
\nonumber
\end{align}
where $K(p)$ is the length of the shortest prefix-free program that computes 
$p(x,y)$ from input $(x,y)$.
\end{lemma}
\begin{remark}\label{rem.sccomp}
Above we required $p(\cdot, \cdot)$ to be computable.
Actually, we only require that $p$ be a
lower semicomputable function, which is a weaker requirement than
recursivity. However, together with the
condition that $p(\cdot, \cdot)$ is a probability distribution,
$\sum_{x,y}  p(x,y) = 1$, this means that $p(\cdot, \cdot)$ is computable,
\cite{LiVi97}, Section 8.1. 
\end{remark}
\begin{proof}
Rewrite the expectation 
\begin{align*}
\sum_x \sum_y p(x,y) I(x:y) \eqa  
\sum_x \sum_y & p(x,y)  [K(x) 
\\& + K(y) - K(x, y)].
\end{align*}
 Define
$\sum_y p(x,y) = p_1 (x)$
and $\sum_x p(x,y) = p_2(y)$
to obtain
\begin{align*}
\sum_x \sum_y p(x,y) I(x:y) \eqa
 \sum_x & p_1 (x) K(x) 
 + \sum_y p_2 (y) K(y) 
\\& - \sum_{x,y} p(x,y) K(x, y).
\end{align*}
Given the program that computes $p$, we can approximate $p_1 (x)$
by a $q_1 (x,y_0) = \sum_{y \leq y_0} p(x,y)$, and
similarly for $p_2$. That is, the
distributions $p_i$ ($i=1,2$) are lower semicomputable, and
by Remark~\ref{rem.sccomp}, therefore, they are computable.
It is known that for every computable probability mass function $q$
we have $H(q) \lea \sum_x q(x) K(x) \lea H(q) + K(q)$, \cite{LiVi97},
Section 8.1. 

Hence, $H(p_i) \lea \sum_x p_i (x) K(x) \lea H(p_i) + K(p_i)$  
($i=1,2$), and $H(p) \lea \sum_{x,y} p (x,y) K(x,y) \lea H(p) + K(p)$.
On the other hand, the probabilistic mutual information
 (\ref{eq.mutinfprob}) is expressed in the entropies by
$I(X;Y) = H(p_1) + H(p_2) - H(p)$.
By construction of the $q_i$'s above,
we have $K(p_1), K(p_2) \lea K(p)$. Since the complexities
are positive, substitution
establishes the lemma.
\end{proof}

Can we get rid of the $K(p)$ error term? The answer is affirmative;
by putting $p(\cdot)$ in the conditional we even get rid of 
the computability requirement.

\begin{lemma}
Given a joint probability mass distribution $p(x,y)$ over $(x,y)$
(not necessarily computable) we have 
\[ I(X;Y)  \eqa  \sum_x \sum_y p(x,y) I(x:y \mid p) , \]
where the auxiliary $p$ means that we can directly access the
values $p(x,y)$ on the
auxiliary conditional information tape of the reference
universal prefix machine.
\end{lemma}

\begin{proof}
The lemma follows from the definition of conditional 
algorithic mutual information, Remark~\ref{rem.cami},
if we show that $\sum_{x}
p(x) K(x \mid p) \eqa H(p)$,
where the $O(1)$ term implicit in the $\eqa$ sign
is independent of $p$.

Equip the reference universal prefix machine,
with an $O(1)$ length
program to compute a Shannon-Fano code from the auxiliary table
of probabilities.
Then, given an input $r$, it can determine
whether $r$ is the Shannon-Fano code word for some $x$.
Such a code word
has length $\eqa - \log p(x)$.
If this is the case, then the machine
outputs $x$, otherwise it halts without output. Therefore,
$K(x   \mid p) \lea - \log p(x)$.
This shows
the upper bound on the expected prefix complexity. 
The lower bound follows as usual
from the Noiseless Coding Theorem.
\end{proof}

We prove a strong version of the information non-increase law
under deterministic processing (later we need the attached corollary):

\begin{theorem}
Given $x$ and $z$, let $q$ be a program 
computing $z$ from $x^*$.
Then
 \begin{equation}\label{eq.nonincrease2}
   I(z : y) \lea I(x : y) + K(q).
 \end{equation}
\end{theorem}

\begin{proof}
By the triangle inequality,
 \begin{align*}
   K(y \mid x^{*}) & \lea K(y \mid z^{*}) + K(z \mid x^{*})
 \\&  \eqa K(y \mid z^{*})+ K(q).
 \end{align*}
Thus,
 \begin{align*}
   I(x : y) & = K(y) - K(y \mid x^{*})
\\ & \gea K(y) - K(y \mid z^{*}) - K(q)
 \\ &  = I(z : y) - K(q).
 \end{align*}
\end{proof}

This also implies the slightly weaker but intuitively
more appealing statement that the mutual information between strings 
$x$ and $y$ cannot be increased by processing $x$ and $y$ separately by
deterministic computations.
 \begin{corollary} Let $f, g$ be recursive functions.
Then
 \begin{equation}\label{eq.nonincrease}
   I(f(x) : g(y)) \lea I(x : y) + K(f)+K(g).
 \end{equation}
 \end{corollary}
 \begin{proof}
It suffices to prove the case $g(y) = y$ and apply it twice.
The proof is by replacing the program $q$ that computes 
a particular string $z$ 
from a particular $x^*$ in (\ref{eq.nonincrease2}). There, $q$
possibly depends on $x^*$ and $z$. Replace it by a program $q_f$ that first
computes $x$ from $x^*$, followed  by computing a
recursive function 
$f$, that is,  $q_f$ is independent of $x$.
Since we only require an $O(1)$-length program to compute
$x$ from $x^*$ we can choose $l(q_f) \eqa K(f)$.

By the triangle inequality,
 \begin{align*}
   K(y \mid x^{*}) & \lea K(y \mid f(x)^{*}) + K(f(x) \mid x^{*})
 \\&  \eqa K(y \mid f(x)^{*})+ K(f).
 \end{align*}
Thus,
 \begin{align*}
   I(x : y) & = K(y) - K(y \mid x^{*}) 
\\ & \gea K(y) - K(y \mid f(x)^{*}) - K(f)
 \\ &  = I(f(x) : y) - K(f).
 \end{align*}
 \end{proof}

It turns out that furthermore, randomized computation can increase
information only with negligible probability.
Let us define the {\em universal probability} $\m(x) = 2^{-K(x)}$.
This function is known to be
maximal within a multiplicative constant among lower semicomputable
semimeasures.
So, in particular, for each computable measure $\nu(x)$ we have 
$\nu(x) \lem \m(x)$, where the constant factor in $\lem$ depends on $\nu$.
This property also holds when we have an extra parameter, like $y^*$, 
in the condition.

Suppose that $z$ is obtained from $x$ by some randomized computation.
The probability $p(z \mid x)$ of obtaining $z$ from $x$ is a semicomputable
distribution over the $z$'s.
Therefore it is upperbounded by 
$\m(z \mid x) \lem \m(z \mid x^{*}) = 2^{-K(z \mid x^{*})}$.
The information increase $I(z : y) - I(x : y)$ satisfies the theorem below.

 \begin{theorem} For all $x,y,z$ we have
 \[
  \m(z \mid x^{*}) 2^{I(z : y) - I(x : y)} 
  \lem \m(z \mid x^{*}, y, K(y \mid x^{*})).
 \]
 \end{theorem}

\begin{remark}
\rm
For example, the probability of an increase of mutual information
by the amount $d$ is $\lem 2^{-d}$.
The theorem 
implies $\sum_{z} \m(z \mid x^{*}) 2^{I(z : y) - I(x : y)} \lem 1$,
the $\m(\cdot \mid x^{*})$-expectation of the exponential of the increase
is bounded by a constant.
\end{remark}

 \begin{proof}
We have
 \begin{align*}
 I(z : y) - I(x : y) & = K(y) - K(y \mid z^{*}) - (K(y) - K(y \mid x^{*}))
 \\&  = K(y \mid x^{*}) - K(y \mid z^{*}).
 \end{align*}
 The negative logarithm of the left-hand side in the theorem is therefore
 \[
  K(z \mid x^{*}) + K(y \mid z^{*}) - K(y \mid x^{*}).
 \]
Using Theorem~\ref{lem.magic}, and the conditional
additivity (\ref{eq.soi-cond}), this is 
 \[
   \gea K(y, z \mid x^{*}) - K(y \mid x^{*}) \eqa 
        K(z \mid x^{*}, y, K(y \mid x^{*})).
 \]
 \end{proof}

\section{Finite Set Models}
\label{s.set}

For convenience, we initially consider the {\em model class} consisting of
the family of finite sets of finite binary strings, that is,
the set of subsets of $\{0,1\}^*$.

\subsection{Finite Set Representations}
\label{ss.exim}
Although all finite sets are recursive there are different
ways to represent or specify the set. 
We only consider ways that have in common
a method of recursively enumerating 
the elements of the finite set one by one, and differ in knowledge
of its size.
For example, we can specify a
set of natural numbers by giving an explicit table or a decision
procedure for membership and a bound on the largest element, 
or by giving a recursive enumeration  
of the elements together with the number of elements,
or by giving a recursive enumeration of the elements
together with a bound on the running time.
We call a representation of a finite set $S$ {\em explicit} 
if the size $|S|$ of the finite set
can be computed from it. 
A representation of $S$ is {\em implicit} if the
logsize $\lfloor \log |S| \rfloor$ can be computed from it.

\begin{example}\label{xmp.implicit}
In Section ~\ref{s.implic}, we will introduce the set 
$S^{k}$ of strings whose elements have complexity $\le k$.
It will be shown that 
this set can be represented implicitly by a program of size $K(k)$,
but can
be represented explicitly only by a program of size $k$.
\end{example}

Such representations are useful in two-stage encodings where 
one stage of the code consists of an index in $S$ of length
$\eqa \log |S|$. 
In the implicit case we know, within an additive constant, how long an
index of
an element in the set is.

We can extend the notion of Kolmogorov complexity from finite
binary strings to finite sets:
The (prefix-) complexity $K_X (S)$ of a
finite set $S$ is defined by
\begin{align*}
K_X(S) = \min_i \{K(i): & \mbox{\rm Turing machine } T_i
\; \; \mbox{\rm computes } S
\\& \mbox{\rm in representation format } X \},
\end{align*}
where $X$ is for example ``implicit'' or ``explicit''.
In general $S^*$ denotes the first shortest self-delimiting binary program
($l(S^*)=K(S)$) in enumeration order
from which $S$ can be computed. These definitions depend,
as explained above, crucial on the representation format $X$: the
way $S$ is supposed to be represented as output of the computation
can make a world of difference for $S^*$ and $K(S)$. 
Since the representation format 
will be clear from the context, and to simplify notation, we
drop the subscript $X$.
To complete our discussion: the worst case of representation format $X$, 
a recursively enumerable representation
where nothing is known about the size of the finite set,
would lead to indices of unknown length.
We do not consider this case.

%pg: 
We may use the notation
 \[
   S_{\impl}, S_{\expl}
 \]
for some implicit and some explicit representation of $S$.
When a result applies to both implicit and explicit representations, 
or when it is clear from the context which representation is meant, we
will omit the subscript.

\subsection{Optimal Model and Sufficient Statistic}
In the following we will distinguish between ``models'' that are 
finite sets, and the ``shortest programs'' to compute those models
that are finite strings. 
Such
a shortest program is in the proper sense a statistic of the data sample
as defined before.
In a way this distinction between ``model'' and ``statistic''
is artificial, but for now we prefer clarity and unambiguousness
in the discussion.

Consider a string $x$
of length $n$ and prefix complexity $K(x)=k$.
We identify the {\em structure} or {\em regularity} in $x$ that are
to be summarized with a set $S$
of which $x$ is a {\em random} or  {\em typical} member: 
given $S$ (or rather, 
%pg: % a
an (implicit or explicit) 
shortest program $S^*$ for $S$), $x$ cannot 
be described significantly shorter than by its maximal length index in $S$, 
that is, $ K(x \mid S^*) \gea \log |S| $. Formally,

\begin{definition}
Let $\beta \ge 0$ be an agreed upon, fixed, constant. 
A finite binary string $x$
is a {\em typical} or {\em random} element of a set $S$ of finite binary
strings if $x \in S$ and 
\begin{equation}\label{eq.deftyp}
 K(x \mid S^*) \ge \log |S| - \beta,
\end{equation}
where $S^*$ is an implicit or explicit shortest program for $S$.
We will not indicate the dependence on $\beta$ explicitly, but the
constants in all our inequalities ($\lea$) will be allowed to be functions
of this $\beta$.
\end{definition}

This definition requires a finite $S$.
In fact, since
$K(x \mid S^*) \lea K(x) $, it limits the size of $S$ to $O(2^k)$
and the shortest program $S^*$ from
which $S$ can be computed) is an {\em algorithmic statistic} for $x$ iff
\begin{equation}\label{eq.typ}
 K(x \mid S^*) \eqa \log |S| .
\end{equation}
Note that the notions of optimality and typicality are not absolute
but depend on fixing the constant implicit in the $\eqa$.
Depending on whether $S^{*}$ is an implicit or explicit program, our
definition splits into implicit and explicit typicality.

%pg:
\begin{example}\label{xmp.typical}
Consider the set $S$ of binary strings of length $n$
whose every odd position is 0.
Let $x$ be an element of this set in which the subsequence of bits in 
even positions is an incompressible string.
Then $S$ is explicitly as well as implicitly typical for $x$.
The set $\{x\}$ also has both these properties.
\end{example}

%pg:
\begin{remark}\label{rem.expl-impl}
It is not clear whether explicit typicality implies implicit typicality.
Section~\ref{s.non-stoch} will show some examples which are implicitly
very non-typical but explicitly at least nearly typical.
\end{remark}

%The index in $S$ is then necessarily random, even conditional to $S$.

There are two natural measures of suitability of such a statistic.
We might prefer either the simplest set, or the largest set, as
corresponding to the most likely structure `explaining' $x$.
The singleton set $\{x\}$, while certainly a statistic for $x$,
would indeed be considered a poor explanation.
Both measures relate to the optimality of a two-stage description of
$x$ using $S$:
\begin{align}\label{eq.twostage}
 K(x) \leq K(x,S)  & \eqa  K(S) + K(x \mid S^*) 
\\ & \lea K(S) + \log |S|,
\nonumber
\end{align}
where we rewrite $K(x,S)$ by~(\ref{eq.soi}).
Here, $S$ can be understood as either $S_{\impl}$ or $S_{\expl}$.
Call a set $S$ (containing $x$)
 for which 
\begin{equation}\label{eq.optim}
K(x) \eqa K(S) + \log |S|,
\end{equation}
{\em optimal}.
%(More precisely, we should require $K(x) \ge K(S) + \log |S| - \beta$.)
Depending on whether $K(S)$ is understood as $K(S_{\impl})$ or
$K(S_{\expl})$, our
definition splits into implicit and explicit optimality.
Mindful of our distinction between a finite set $S$ and a
program that describes $S$ in a required representation format,
we call a shortest program for an optimal set with respect to $x$
an {\em algorithmic sufficient statistic} for $x$.
Furthermore, among optimal sets,
there is a direct trade-off between complexity and logsize, which together
sum to $\eqa k$. 
Equality (\ref{eq.optim}) is the algorithmic equivalent
dealing with the relation between the individual
sufficient statistic and the individual data sample, in contrast
to the probabilistic notion (\ref{eq.suffstatprob}).

\begin{example}
\label{ex.restricted}
The following restricted model family illustrates the difference
between the algorithmic individual notion of sufficient
statistic and the probabilistic averaging one.
Foreshadowing the discussion in section \ref{s.discussion}, 
this example also illustrates the idea that the semantics of the model class
should be obtained by a restriction on the family of allowable
models, after which the (minimal) sufficient statistic identifies
the most appropriate model in the allowable family and thus optimizes
the parameters in the selected model class.
In the algorithmic setting we use all subsets of $\{0,1\}^n$
as models and the shortest programs computing them from a
given data sample as the statistic. Suppose we have background
information constraining the family of models to the $n+1$
finite sets $S_s = \{x \in \{0,1\}^n : x= x_1 \ldots x_n \& 
\sum_{i=1}^n x_i = s \}$ ($0 \leq s \leq n$).
Assume that our model family is the family of Bernoulli distributions.
Then, in the
probabilistic sense for every data sample $x = x_1 \ldots x_n$ there
is only one natural sufficient statistic: for $\sum_i x_i = s$
this is $T(x)=s$ with the corresponding model $S_s$.
In the algorithmic setting the situation is more subtle. (In
the following example we use the complexities conditional on $n$.)
For $x=x_1 \ldots x_n$ with $\sum_i x_i = \frac{n}{2}$
taking $S_{\frac{n}{2}}$ as model yields 
$|S_{\frac{n}{2}}| = {n \choose {\frac{n}{2}}}$, 
and therefore $\log | S_{\frac{n}{2}} | \eqa n - \frac{1}{2} \log n$.
The sum of $K(  S_{\frac{n}{2}} |n) \eqa 0$ and the logarithmic
term  gives $ \eqa n - \frac{1}{2} \log n$ for the right-hand side of
(\ref{eq.optim}). But taking $x = 1010 \ldots 10$ yields $K(x \mid n) \eqa 0$
for the left-hand side. Thus, there is no algorithmic sufficient 
statistic for the latter $x$ in this model class, while every $x$
of length $n$ has a probabilistic sufficient statistic in the model class.
In fact, the restricted model class has algorithmic sufficient
statistic for data samples $x$ of length $n$ that have maximal complexity
with respect to the frequency of ``1''s, the other data samples
have no algorithmic sufficient statistic in this model class.
\end{example}

%pg:
\begin{example}\label{xmp.optimal}
It can be shown that the set $S$ of Example~\ref{xmp.typical} is also
optimal, and so is $\{x\}$.
%pg: I shifted the position of this example
Typical sets form a much wider class than optimal ones:
$\{x,y\}$ is still typical for
$x$ but with most $y$, it will be too complex to be optimal for $x$.

For a perhaps less artificial example, consider complexities conditional
on the length $n$ of strings.
Let $y$ be a random string of length $n$, let
$S_{y}$ be the set of strings of length $n$ which have 0's exactly
where $y$ has, and let $x$ be a random element of $S_{y}$.
Then $x$ is a string random with respect to the distribution in which 
1's are chosen independently 
with probability 0.25, so its complexity is much less than $n$.
The set $S_{y}$ is typical with respect to $x$ but is too complex to be optimal,
since its (explicit or implicit) complexity conditional on $n$ is $n$.
\end{example}

It follows that (programs for) optimal sets are statistics.
Equality (\ref{eq.optim}) expresses the conditions on the algorithmic individual
relation between the data and the sufficient statistic. Later we 
demonstrate that this relation implies that the probabilistic
optimality of mutual information (\ref{eq.mutinfprob}) holds
for the algorithmic version in the expected sense.

An algorithmic sufficient statistic $T( \cdot)$
is a sharper individual notion than a probabilistic sufficient
statistic. An optimal set $S$ associated with $x$ (the shortest
program computing $S$ is the corresponding
sufficient statistic associated with $x$) is chosen such that
$x$ is maximally random with respect to it. That is, the 
information in $x$ is divided in a relevant structure expressed
by the set $S$, and the remaining randomness with respect
to that structure, expressed by $x$'s index in $S$ of $\log |S|$
bits. The shortest program for $S$ is itself alone an algorithmic
definition of structure, without a probabilistic interpretation.
%%For every
%%joint probability mass distribution the 
%%expectation of the algorithmic mutual information $I(\theta : x)$
%%is close to the expectation of $I(\theta : T(x))$,
%%corresponding to the probabilistic sufficient
%%statistic condition (\ref{eq.suffstatprob}).
%%If this were not the case then the algorithmic notion would not
%%be a sharpening of the probabilistic notion to individual objects,
%%but something else. For convenience we 
%%let $\theta$ be the parameter (vector) taking values
%%in the set of natural numbers ${\cal N}$, and let $\Omega = \{0,1\}^*$
%%be the sample space.
%%The equality holds in the following sense:
%%
%%\begin{lemma}
%%Given an algorithmic sufficient 
%%statistic $T: \{0,1\}^* \rightarrow 2^{\Omega} $
%%and a computable joint probability mass 
%%distribution $p(\theta,x)$ over ${\cal N} \times \Omega$
%%we have $ \sum_{\theta} \sum_x p(\theta,x) I(\theta:x) \eqa 
 %%\sum_\theta \sum_x p(\theta,x) I(\theta : T(x))$ up to an
%%additive term linear in $K(p)$,
%%%\[ I(x;y) \lea  \sum_x \sum_y p(x,y) I(x:y) \lea I(x;y) + 3 K(p) , \]
%%where $K(p)$ is the length of the shortest prefix-free program that computes
%%$p(\theta ,x)$ from input $(\theta ,x)$.
%%\end{lemma}
%%\begin{proof}
%%%%By (\ref{def.mutinf}) we only need to show 
 %%%%$ \sum_{\theta} \sum_x p(\theta,x) K(\theta \mid x^*) \eqa
 %%\sum_\theta \sum_x p(\theta,x) K(\theta \mid T(x)^*)$.
%%Since $K(\theta |x^*) \lea K(\theta |T(x)^*)$ ???
%%\end{proof}

One can also consider notions of 
{\em near}-typical and {\em near}-optimal that arise from replacing
the $\beta$  in (\ref{eq.deftyp})
by some slowly growing functions, such as $O(\log l(x))$ or
$O(\log k)$ as in~\cite{ShenStoch83,ShenStoch99}.

In~\cite{ShenStoch83,Vy87}, a function 
of $k$ and $x$ is defined as the lack of typicality
of $x$ in sets of complexity at most $k$, and they then consider the
minimum
$k$ for which this function becomes $\eqa 0$ or very small. This is
equivalent to our notion of a typical set.
%pg:
See the discussion of this function in Section~\ref{s.non-stoch}.
In~\cite{Co85,CT91}, only optimal sets are considered, and the one
with the shortest program
is identified as the {\em algorithmic minimal sufficient statistic} of $x$.
Formally, this is the shortest program that computes
a finite set $S$ such that~(\ref{eq.optim}) holds.

%Summarizing: all algorithmic minimal sufficient 
%statistics are sufficient statistics,
%and all sufficient statistics are typical statistics. 

\subsection{Properties of Sufficient Statistic}

We start with a sequence of lemmas that will be used in the later
theorems.
Several of these lemmas have two versions: for implicit sets and for explicit
sets.
In these cases, $S$ will denote $S_{\impl}$ or $S_{\expl}$ respectively.

Below it is shown that the mutual information between every 
typical set and the
data is not much less than $K(K(x))$, the complexity of the complexity
$K(x)$ of the data $x$. 
For optimal sets it is at least that, and for
algorithmic minimal statistic it is equal to that.
%pg: I do not understand the following sentence.
% This suggests that the mutual information must be larger
% for more ``random'' data complexities.
The number of elements
of a typical set is determined by the following:

\begin{lemma}\label{lem.card}
Let $k=K(x)$.
If a set $S$ is
(implicitly or explicitly)
typical for $x$ then 
$I(x:S) \eqa k- \log |S|$.
%$|S| = \Theta (2^{k-I(x:S)})$.
\end{lemma}
\begin{proof}
By definition $I(x:S) \eqa K(x) - K(x \mid S^*)$
and by typicality $K(x \mid S^*) \eqa \log |S|$.
\end{proof}

Typicality, optimality, and minimal optimality
successively restrict the range of the cardinality (and complexity)
of a corresponding model for a data $x$.
The above  lemma states that for 
(implicitly or explicitly)
typical $S$ the cardinality $|S| = \Theta (2^{k-I(x:S)})$. 
The next lemma asserts that for implicitly typical $S$
the value $I(x:S)$ can fall below $K(k)$ by no more than an
additive logarithmic term.

\begin{lemma}\label{lem.KI}
Let $k=K(x)$.
If a set $S$ is (implicitly or explicitly)
typical for $x$ then $I(x:S) \gea K(k) - K(I(x:S))$
and $\log |S| \lea k-  K(k) + K(I(x:S))$.
(Here, $S$ is understood as $S_{\impl}$ or $S_{\expl}$ respectively.)
\end{lemma}
\begin{proof}
Writing $k=K(x)$, since 
\begin{equation}\label{eq.kx}
k \eqa K(k,x) \eqa K(k)+K(x \mid k^*)
\end{equation}
 by~(\ref{eq.soi}), we have
$I(x:S) \eqa K(x)-K(x \mid S^*) \eqa K(k)-[K(x \mid S^*)-K(x \mid k^*)]$.
Hence, it suffices to show $K(x \mid S^*)-K(x \mid k^*) \lea K(I(x:S))$.
Now, from an implicit description $S^*$ we can
find the value $\eqa \log |S|\eqa k-I(x:S)$. To recover $k$ we only
require an extra $K(I(x:S))$ bits apart from $S^*$.
Therefore, $K(k \mid S^*) \lea K(I(x:S))$.
This reduces what we have to show to
$K(x \mid S^*) \lea K(x \mid k^*) + K(k \mid S^*)$ which is asserted by
Theorem~\ref{lem.magic}.

%%Since $K(I(x:S))$ cannot be negative,
%%if $I(x:S) \gea K(k)$ then there is nothing to prove. Hence
%%we assume $I(x:S) \lea K(k)$. 

%By definition $I(x:S) \eqa k-K(x|S^*)$ and substituting
%the above expressions in the right hand side we find
%$I(x:S)$ by typicality
%$K(x|S^*) \eqa \log |S|$. Therefore, 
%%By Lemma~\ref{lem.card}, $\log |S| + I(x:S) \eqa k$.
 %%Clearly  $K(x) = K(x,K(x))$, and therefore
%%$k \eqa  K(k)+ K(x \mid k^*)$ by~(\ref{eq.soi}).  
%%
%%Moreover,
%%$K(x \mid S^*) \lea K(x,k \mid S^*) \lea K(x \mid k, K(k \mid S^*),S^*) + K(k \mid S^*) 
%%\lea K(x \mid k, K(k \mid S^*),S^*)+K(k \mid S^*)
%%\lea K(I(x:S)) + K(k \mid S^*)$. The last inequality holds since
%%from an implicit description $S^*$ we can 
%%find $\eqa \log |S|\eqa k-I(x:S)$; hence to recover $k$ we only
%%require an extra $K(I(x:S))$ bits apart from $S^*$.
%%
%%Substituting these bounds
%%on $k$ and $K(x \mid S^*)$ in $I(x:S) \eqa k-K(x \mid S^*)$ gives 
%%$I(x:S) \gea K(k) - K(I(x:S))$. By Lemma~\ref{lem.card}
%%we then find the upper bound on $|S|$.
\end{proof}

The term $I(x:S)$ is at least $K(k)- 2 \log K(k)$ where $k=K(x)$. 
For $x$ of length $n$ with
$k \gea n$ and $K(k) \gea l(k) \gea \log n$,
this yields $I(x:S) \gea \log n - 2 \log \log n$.
%pg: This last sentence is not an interesting observation.

If we further restrict typical sets to optimal sets then
the possible number of elements in $S$ is slightly restricted.
First  we show that implicit optimality of a set with respect to a data
 is equivalent to typicality with respect to the data
combined with effective constructability (determination) from the data.

\begin{lemma}\label{lem.ioty}
A set $S$ is 
(implicitly or explicitly)
optimal for $x$ iff it is typical and $K(S \mid x^*) \eqa 0$.
\end{lemma}

\begin{proof}
A set $S$ is optimal iff~(\ref{eq.twostage}) holds with equalities.
Rewriting $K(x,S) \eqa K(x)+K(S \mid x^*)$
the first inequality becomes an equality iff $K(S \mid x^*) \eqa 0$,
and the second inequality becomes an equality iff $K(x \mid S^*) \eqa \log |S|$
(that is, $S$ is a typical set).
\end{proof}

\begin{lemma}\label{lem.opt}
Let $k=K(x)$.
If a set $S$ is (implicitly or explicitly) optimal for $x$, then 
$I(x:S) \eqa K(S) \gea K(k)$ and $\log |S| \lea  k-K(k)$.
\end{lemma}

\begin{proof}
If $S$ is optimal for $x$, then 
$k = K(x) \eqa K(S)+K(x \mid S^*) \eqa K(S)+\log |S|$.
%It follows that $K(k \mid S^*) \eqa K(K(S)+(\log |S|) \mid S^*) \eqa 0$,
%hence $K(k) \lea K(S) + K(k \mid S^*) \eqa K(S)$. 
  From $S^*$ we can find both $K(S) \eqa l(S^*)$
and $\eqa \log |S|$ and hence $k$, that is, $K(k) \lea K(S)$.
We have $I(x:S) \eqa K(S)-K(S \mid x^*) \eqa K(S)$ by~(\ref{eq.soi}),
Lemma~\ref{lem.ioty}, respectively.
This proves the first property.
Substitution of $I(x:S) \gea K(k) $ in the expression
of Lemma~\ref{lem.card} proves the second property.
\end{proof}

\begin{figure}
\begin{center}
\epsfxsize=8cm
\epsfxsize=8cm \epsfbox{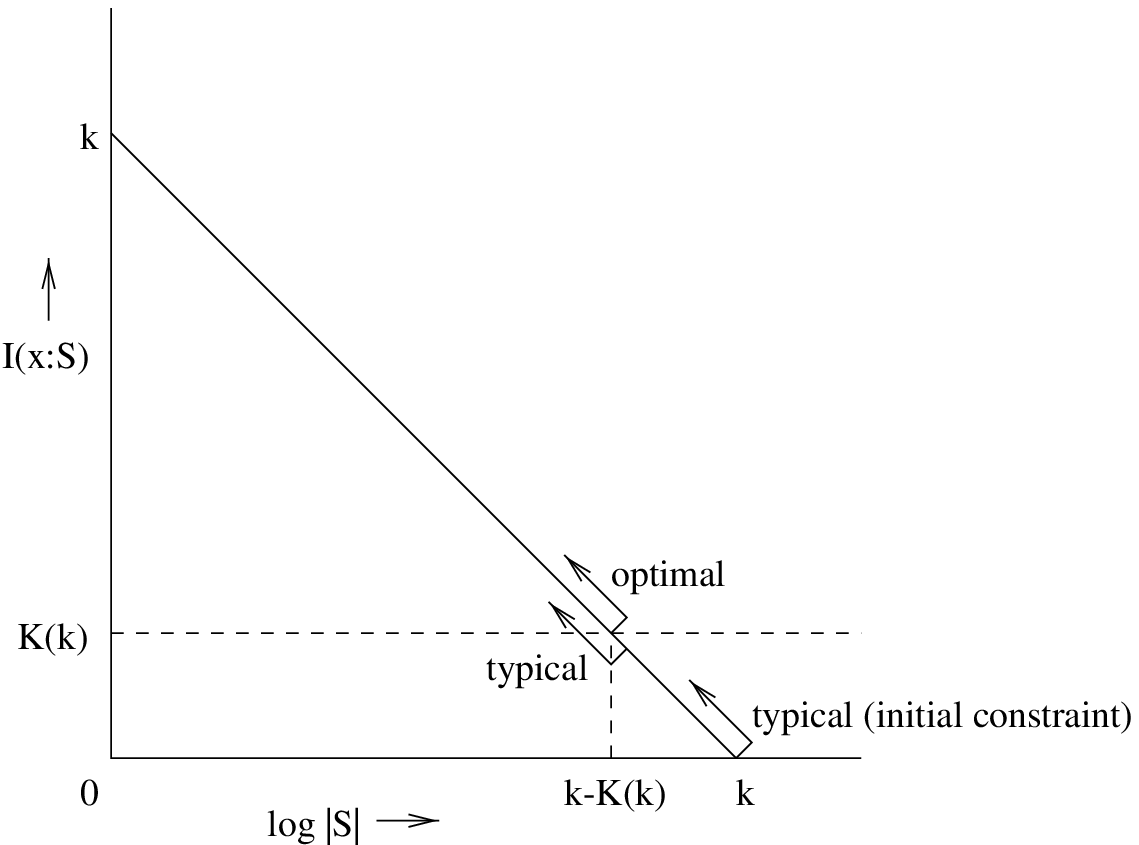}
\end{center}
\caption{Range of statistic on the straight line $I(x:S)
\eqa K(x)-\log |S|$.}
\label{stat}
\end{figure}

\subsection{Implicit Minimal Sufficient Statistic}\label{s.implic}

A simplest implicitly optimal set (that is, of least complexity)
is an implicit algorithmic minimal sufficient statistic. We demonstrate that
$S^k = \{y: K(y) \leq k\}$, the set of all
strings of complexity at most $k$, is such a set. First we establish the
cardinality of $S^k$:

%pg: I think this is proven by Solovay (see my lecture notes
% and probably your book)
\begin{lemma}\label{lem.impl}
$\log |S^k| \eqa k-K(k)$.
\end{lemma}

\begin{proof}
The lower bound is easiest.
Denote by $k^*$ of length $K(k)$
a shortest program for $k$. Every string $s$ of length
$k-K(k)-c$ can be described in a self-delimiting manner
by prefixing it with $k^*c^*$, hence $K(s) \lea k-c+2\log c $.
%By considering strings $s$ of length $k-K(k)-c$ instead, 
For a large
enough constant $c$, we have $K(s) \leq k$
and hence there are $\Omega(2^{k-K(k)})$ strings that are in $S^k$.

For the upper bound: by~(\ref{eq.kx}),
all $x \in S^k$ satisfy $K(x \mid k^*) \lea k-K(k)$,
and there can only be $O(2^{k-K(k)})$ of them.
\end{proof}

   From the definition of $S^k$ it follows that it is defined by $k$
alone, and it is the same set that is optimal for all objects
of the same complexity $k$.

\begin{theorem}
The set $S^k$ is implicitly optimal for every $x$ with $K(x)=k$.
%pg:
Also, we have $K(S^{k}) \eqa K(k)$.
\end{theorem}
\begin{proof}
%Fix a constant $c$ such that $|S^k| \leq 2^{k-K(k)+c}$.
%A string $I^k_x$ of length $\eqa k-K(k)$ denotes
%the index of $x$ in the standard enumeration of $S^k$
%(Lemma~\ref{lem.impl}). Then $k^* I^k_x$ together with
%a fixed decoding program is a self-delimiting description of $x$
%of length $\eqa k$. This shows $K(x) \eqa K(k) + \log |S^k|$.
  From $k^*$ we can compute both $k$  and $k-l(k^*) = k-K(k)$
and recursively enumerate $S^k$.
Since also $\log |S^k| \eqa k-K(k)$ (Lemma~\ref{lem.impl}),
the string 
$k^*$ plus a fixed program
is an implicit description of $S^k$ so that $K(k) \gea K(S^k)$.
Hence, $K(x) \gea K(S^k)+ \log |S^k|$ and since $K(x)$
is the shortest description by definition equality ($\eqa$) holds. 
That is, $S^k$ is optimal for $x$.
By Lemma~\ref{lem.opt}
$K(S^k) \gea K(k)$ which together with the reverse inequality above
yields $K(S^k)  \eqa K(k)$ which shows the theorem.
\end{proof}

Again using Lemma~\ref{lem.opt} shows that the optimal set $S^k$
has least complexity among all optimal sets for $x$, and therefore:

\begin{corollary}\label{corr.iamss}
The set $S^k$ is an implicit algorithmic minimal sufficient statistic
for every $x$ with $K(x)=k$.
\end{corollary}

All algorithmic minimal sufficient statistics 
$S$ for $x$ have $K(S) \eqa K(k)$, and therefore there are
$O(2^{K(k)})$ of them. At least one such a statistic ($S^k$)
is associated with every one of the $O(2^k)$ strings $x$ of
complexity $k$.
Thus, while the idea of the algorithmic minimal sufficient
statistic is intuitively appealing, its unrestricted
use doesn't seem to uncover most relevant aspects of reality.
The only relevant structure in the data
with respect to an algorithmic minimal sufficient statistic 
is the Kolmogorov complexity.
To give an example, an initial 
segment of $3.1415 \ldots$ of length $n$ of complexity
$\log n + O(1)$ shares the same algorithmic sufficient statistic
with many (most?) binary strings of length $\log n + O(1)$.

\subsection{Explicit Minimal Sufficient Statistic}
Let us now consider representations of finite sets that are explicit
in the sense that we can compute the cardinality of the set from the
representation.

\subsubsection{Explicit Minimal Sufficient Statistic: Particular Cases}
\begin{example}
The description program
enumerates all the elements of the set and halts.
Then a set like
$S^k  = \{y : K(y) \leq k \}$ has complexity $\eqa k$~\cite{ShenStoch99}:  
Given the program we can find an element 
not in $S^k$, which element by definition has complexity $>k$.
Given $S^k$ we can find this element
and hence $S^k$ has complexity $\gea k$.
Let 
 \[
  N^k = |S^k|,
 \]
then by Lemma~\ref{lem.impl} $\log N^{k} \eqa k-K(k)$.
We can list $S^k$ given $k^{*}$ and $N^{k}$ 
which shows $K(S^k) \lea k$.
\end{example}

\begin{example}
One way of implementing explicit finite representations is to
provide an explicit generation time for the enumeration process.
If we can generate $S^k$ in time $t$ recursively
using $k$, then the previous argument shows that the complexity
of every number $t' \geq t$ satisfies $K(t',k) \geq k$ so that
$K(t') \gea K(t' \mid k^*) \gea k - K(k)$ by~(\ref{eq.soi}). 
This means that 
$t$ is a huge time which as a function of $k$ rises faster
than every computable function.
This argument also shows
that explicit enumerative descriptions of sets $S$ containing $x$ by
an enumerative process $p$ plus a limit on the computation time $t$
may take only $l(p)+ K(t)$ bits 
(with $K(t) \leq \log t + 2 \log \log t$)
but $\log t$ unfortunately becomes {\em noncomputably large}!
\end{example}

\begin{example}
Another way is to indicate the element of $S^k$ that requires
the longest generation time as part of the dovetailing process,
for example by its index $i$ in the enumeration, $i \leq 2^{k-K(k)}$.
Then, $K(i \mid k) \lea k - K(k)$.
In fact, since a shortest program $p$ for the $i$th element
together with $k$ allows us to generate $S^k$ explicitly,
and abive we have seen that explicit description format
yoelds $K(S^k) \eqa k$, we find
we have $K(p,k) \gea k$ and hence $K(p) \gea k-K(k)$.
\end{example}

In other cases the generation time is simply recursive in the input:
$S_n = \{y : l(y) \leq n \}$ so that 
$K(S_n) \eqa K(n) \leq \log n + 2 \log \log n$. 
That is, this sufficient
statistic for a random string $x$ with 
$K(x)\eqa n+K(n)$ has complexity $K(n)$ both for implicit descriptions and
explicit descriptions: differences in complexity arise 
only for nonrandom strings (but not too nonrandom, for
$K(x) \eqa 0$ these differences vanish again). 
\begin{lemma}
$S_n$ is an example of a 
{\em minimal sufficient statistic}, both explicit and implicit,
for all $x$ with $K(x) \eqa n + K(n)$.
\end{lemma}
\begin{proof}
The set $S_n$ is a sufficient statistic for $x$ since
$K(x) \eqa K(S_n)+ \log |S_n|$. It is minimal since
by Lemma~\ref{lem.opt} we must have $K(S) \gea K(K(x))$ for
implicit, and hence for
explicit sufficient statistics. It is evident that $S_n$ is explicit:
$|S_n|=2^n$.
\end{proof}

It turns out that some strings cannot thus be explicitly
represented parsimonously
with low-complexity models 
(so that one necessarily has bad high complexity
models like $S^k$ above).
For explicit representations, 
\cite{ShenStoch83}
has demonstrated the
existence of a class of strings called {\em non-stochastic}
that don't have efficient two-part representations
with $K(x) \eqa K(S) + \log |S|$   ($x \in S$) with $K(S)$
significantly less than  $K(x)$. This result does not yet
enable us to exhibit an explicit minimal sufficient statistic
for such a string. 
But in Section~\ref{s.non-stoch} we improve
these results to the best possible, simultaneously establishing
explicit minimal sufficient statistics for the subject
ultimate non-stochastic strings:
\begin{lemma}\label{lem.emss}
For every length $n$, there exist strings $x$ of length $n$
with $K(x \mid n) \eqa n$ for which $\{x\}$ is an explicit
minimal sufficient statistic.
The proof 
is deferred to the end of Section~\ref{s.non-stoch}.
\end{lemma}

\subsubsection{Explicit Minimal Near-Sufficient Statistic: General Case}
Again, consider the special set $S^k = \{y : K(y) \leq k \}$.
As we have seen earlier, $S^k$ itself cannot be explicitly
optimal for $x$ since $K(S^k) \eqa k$ and $\log N^{k} \eqa k-K(k)$,
and therefore $K(S^k)+\log N^{k} \eqa 2k-K(k)$ which considerably
exceeds $ k$.
However, it turns out that a closely related set ($S^k_{m_x}$
below) is explicitly near-optimal.
Let $I^k_y$ denote the index of $y$ in the standard enumeration of $S^k$,
where all indexes are padded to the same length $\eqa k-K(k)$ with 0's in
front.
For $K(x) = k$, let $m_x$ denote the longest joint prefix of
$I^k_x$ and $N^k$, and let
 \[
  I^k_x = m_x 0 i_x, \quad N^k = m_x 1 n_x.
 \]

\begin{lemma}\label{lem.nearub}
For $K(x) = k$, the set 
$S^k_{m_x} = \{y\in S^k: m_x 0 \; \;\mbox{\rm  a prefix of } I^k_y \}$
satisfies
 \begin{align*}
   \log |S^k_{m_x}| &\eqa k-K(k) - l(m_x),
\\     K(S^k_{m_x}) &\lea K(k) + K(m_x) \lea K(k)+l(m_x)+K(l(m_x)).
 \end{align*}
Hence it is explicitly near-optimal 
for $x$ (up to an addive $K(l(m_x)) \lea K(k) \lea \log k + 2 \log \log k$
term).
\end{lemma}

\begin{proof}
%%If $l(m_x)=0$ then $m_x0i_x = 0 \ldots$ while $N^k=m_x1n_x = 1 \ldots$
%%(since the index of $x$ cannot exceed the count $N^k$). In this case
%%the set $S^k_{m_x} = \{y\in S^k : I_y = 0 \ldots \}$. The remainder
%%of the argument is as in the case $l(m_x) >0$.
%%
%%Assume $l(m_x) >0$.
We can describe $x$ by $k^*m_x^*  i_x$ where $m_x 0 i_x$ is the
index of $x$ in the enumeration of $S^k$. 
%%This total description
%%has length $\eqa k$ since $l(m_x) \eqa 0$ by
%%Lemma~\ref{lem.mx}. 
Moreover, $k^* m_x^*$ explicitly 
describes the set $S^k_{m_x}$. Namely, using $k$ we can recursively
enumerate $S^k$.
At some point the first string $z \in S^k_{m_x}$ is enumerated (index 
$I_z^k = m_x 00 \ldots 0$). 
By assumption $I_x^k = m_x 0 \ldots$ and
$N^k = m_x 1 \ldots $.
Therefore, in the enumeration of $S^k$
eventually string $u$ with $I_u^k = m_x 011 \ldots 1$ occurs
which is the last string in the enumeration of $S^k_{m_x}$.
Thus, the size
of $S^k_{m_x}$ is precisely $2^{l(N^k)-l(m_x)}$, where
$l(N^k)-l(m_x) \eqa l(n_x) \eqa \log |S^k_{m_x}|$, and $S^k_{m_x}$
is explicitly described by $k^*m_x^*$.
Since $l(k^*m_x0 i_x) \eqa k$
and $\log |S^k_{m_x}| \eqa k-K(k)-l(m_x)$ we have 
 \begin{align*}
  K(S^k_{m_x})+ \log |S^k_{m_x}| &\eqa K(k) + K(m_x) +k -K(k)-l(m_x)
\\           &\eqa k+ K(m_x)-l(m_x) \lea k+K(l(m_x)).
 \end{align*}
This shows $S^k_{m_x}$ is explicitly near optimal for $x$ (up
to an additive logarithmic term).
\end{proof}

%%The explicitly optimal set $S^k_{m_x}$ has complexity $\eqa K(k)$
%%(since $K(m_x) \eqa 0$ by Lemma~\ref{lem.mx}).
%%Again using Lemma~\ref{lem.opt} shows
%%that it 
%%has least complexity $\gea K(k)$ among all
%%(even implicitly) optimal sets for $x$, and therefore:
\begin{lemma}\label{lem.nearlb}
Every explicit optimal set $S \subseteq S^k$
containing $x$ satisfies 
 \[
  K(S) \gea K(k) + l(m_x) - K(l(m_x)).
 \]
\end{lemma}
\begin{proof}
If $S \subseteq S^k$ is explicitly optimal for $x$, then we can
find $k$ from $S^*$ (as in the proof of Lemma~\ref{lem.opt}),
and given $k$ and $S$ we find $K(k)$ as in Theorem~\ref{lem.magic}.
Hence, given $S^*$, we can enumerate $S^k$ and determine the
maximal index $I^k_y$ of a $y \in S $.
Since also $x \in S $, the numbers
$I^k_y, I^k_x, N^k$ have a maximal common prefix $m_x$.
Write
$I^k_x = m_x 0 i_x$ with $l(i_x) \eqa k-K(k)-l(m_x)$ by
Lemma~\ref{lem.opt}.
Given $l(m_x)$ we can determine $m_x$ from $I^k_y$. Hence, from
$S,l(m_x)$, and $i_x$ we can reconstruct $x$. 
That is,
$K(S)+K(l(m_x)) + l(I^k_x)- l(m_x) \gea k$, which yields the lemma.
\end{proof}

Lemmas~\ref{lem.nearub},~\ref{lem.nearlb} demonstrate:

\begin{theorem}\label{theo.ekmss}
The set $S^k_{m_x}$ is an explicit 
algorithmic minimal near-sufficient statistic
for $x$ among subsets of $S^k$ in the following sense:
 \begin{align*}
   |K(S^k_{m_x}) - K(k) - l(m_x) | &\lea K(l(m_x)),
\\                \log |S^k_{m_x}| &\eqa k - K(k) - l(m_x).
 \end{align*}
Hence $K(S^k_{m_x})+\log |S^k_{m_x}| \eqa k \pm K(l(m_x))$.
Note, $K(l(m_x)) \lea \log k + 2 \log \log k$.
\end{theorem}

\subsubsection{Almost Always ``Sufficient''}
We have not completely succeeded in giving a concrete algorithmic
explicit minimal sufficient statistic. However, we can show
that $S^k_{m_x}$ is {\em almost always} minimal sufficient.

The complexity and cardinality of $S^k_{m_x}$ depend on $l(m_x)$
which will in turn depend on $x$.
One extreme is $l(m_x) \eqa 0$ which happens for the majority of $x$'s 
with $K(x) = k$---for example, the first 99.9\% in the enumeration order. 
For those $x$'s we can replace ``near-sufficient'' by 
``sufficient'' in Theorem~\ref{theo.ekmss}.
Can the other extreme be reached?
%Clearly, $x$ cannot be enumerated 
%close to the start of the enumeration process of $S^k$ since then
%the complexity of $x$ would drop below $k$. 
This is the case when
$x$ is enumerated close to the end of the enumeration
of $S^k$.
For example, this happens for the ``non-stochastic'' objects
of which the existence was proven by Shen~\cite{ShenStoch83} (see
Section~\ref{s.non-stoch}).
For such objects, $l(m_x)$ grows to $\eqa k-K(k)$ and
the complexity of $S^k_{m_x}$ rises to $\eqa k$
while $\log |S^k_{m_x}|$ drops to $\eqa 0$.
That is, the explicit
algorithmic minimal sufficient statistic for $x$ is essentially 
$x$ itself.
For those $x$'s we can also replace ``near-sufficient'' with ``sufficient''
in Theorem~\ref{theo.ekmss}.
Generally:
%pg: This is too imprecise to be called a formal corollary.
% \begin{corollary}
for the overwhelming majority of data $x$ of complexity $k$
the set $S^k_{m_x}$ is an explicit algorithmic minimal sufficient statistic
among subsets of $S^k$ (since $l(m_x) \eqa 0$).
% \end{corollary}

The following discussion will put what was said above into a
more illuminating context.
Let 
 \[
   X(r) = \setof{x : l(m_{x}) \ge r}.
% Gacs: changed from \geq
 \]
The set $X(r)$ is infinite, but we can break it into slices and bound each
slice separately.

 \begin{lemma}
 \[
   |X(r) \bigcap (S^{k} \setminus S^{k-1}) | \leq 2^{-r+1} |S^{k}|.
% Gacs: changed from \lem 2^{-r} |S^{k}|.
 \]
 \end{lemma}
 \begin{proof}

For every $x$ in the set defined by the left-hand side of the inequality,
we have $l(m_{x}) \ge r$, and the
length of continuation of $m_{x}$ to the total padded index of $x$ is 
$\le \cei{\log |S^{k}|} - r \le \log |S^{k}| - r  + 1$.
% Gacs: changed from just $\leq \log |S^{k}| - r$.
Moreover, all these indices share the same first $r$ bits.
This proves the lemma.
 \end{proof}

 \begin{theorem}\label{t.X(r)-bd}
 \[
   \sum_{x \in X(r)} 2^{-K(x)} \leq 2^{-r+2}.
% Gacs: changed from  \sum_{x \in X(r)} 2^{-K(x)} \lem 2^{-r}.
% Rewrote the proof to prove the more precise statement.
 \]
 \end{theorem}
 \begin{proof}
  Let us prove first
 \begin{equation}\label{e.Sk.bd}
 \sum_{k \ge 0} 2^{-k} |S^{k}| \le 2.
 \end{equation}
By the Kraft inequality, we have, with $t_{k} = |S^{k} \setminus S^{k-1}|$,
 \begin{equation*}
 \sum_{k \ge 0} 2^{-k} t_{k} \le 1,
 \end{equation*}
since $S^{k}$ is in 1-1 correspondence with the prefix programs of length
$\le k$.
Hence
 \[
   \begin{split}
  \sum_{k \ge 0} 2^{-k}|S^{k}| &= \sum_{k \ge 0} 2^{-k} \sum_{i = 0}^{k} t_{i}
      = \sum_{i \ge 0} t_{i} \sum_{k = i}^{\infty} 2^{-k}
\\                             &= \sum_{i \ge 0} t_{i} 2^{-i+1} \le 2.
   \end{split}
 \]
For the statement of the lemma, we have
 \begin{align*}
   \sum_{x \in X(r)} 2^{-K(x)} &=
  \sum_{k \ge 0} 2^{-k} |X(r) \bigcap (S^{k} \setminus S^{k-1})|
\\  &\leq 2^{-r+1} \sum_{k \ge 0} 2^{-k}|S^{k}| \le 2^{-r+2},
 \end{align*}
where in the last inequality we used~(\ref{e.Sk.bd}).
 \end{proof}

This theorem can be interpreted as follows, (we rely here on a discussion,
unconnected with the present topic, about universal probability
with L.~A.~Levin in 1973).
The above theorem states $\sum_{x \in X(r)} \m(x) \leq 2^{-r+2}$.
% Gacs: changed from \lem 2^{-r}.
By the multiplicative dominating property of $\m(x)$ with respect
to every lower semicomputable semimeasure,
it follows that for every computable measure $\nu$, we have 
$\sum_{x \in X(r)} \nu(x) \lem 2^{-r}$.
Thus, the set of objects $x$ for which $l(m_{x})$ is large
has small probability with
respect to every computable probability distribution.

To shed light on the exceptional nature of strings $x$ with large $l(m_{x})$
from yet another direction, let $\chi$ be the infinite binary sequence,
the {\em halting sequence},
which constitutes the characteristic function of 
the halting problem for our universal Turing machine: the $i$th bit
of $\chi$ is 1 of the machine halts on the $i$th program, and is 0 otherwise.
The expression 
 \[
   I(\chi : x) = K(x) - K(x \mid \chi)
 \]
shows the amount of information in the halting sequence about the string
$x$.
(For an infinite sequence $\eta$, we go back formally to the definition
$I(\eta : x) = K(x) - K(x \mid \eta)$ of~\cite{LiVi97}, since
introducing a notion of $\eta^{*}$ in place of $\eta$ here
has not been shown yet to bring any benefits.)
We have
 \[
   \sum_{x} \m(x) 2^{I(\chi : x)} = \sum_{x} 2^{-K(x \mid \chi)} \le 1.
 \]
Therefore, if we introduce a new quantity $X'(r)$ related to $X(r)$ 
defined by
 \[
   X'(r) = \setof{x : I(\chi : x) > r},
 \]
then by Markov's inequality,
 \[
   \sum_{x \in X'(r)} \m(x) 2^{I(\chi : x)} < 2^{-r}.
 \]
That is, the universal probability of $X'(r)$ is small.
This is a new reason for $X(r)$ to be small,
as is shown in the following theorem.

 \begin{theorem}
We have
 \[
   I(\chi : x) \gea l(m_{x}) - 2\log l(m_{x}),
 \]
and (essentially equivalently) $X(r) \subset X'(r - 2\log r)$.
 \end{theorem}

\begin{remark}
The first item in the theorem implies: If $l(m_x) \geq r$,
then $I(\chi : x) \gea r - 2 \log r$. This in its turn implies
the second item $X(r) \subset X'(r - 2\log r)$. Similarly, the second
item essentially implies the first item.
Thus, a string for which the explicit 
minimal sufficient statistic has complexity
much larger than $K(k)$ (that is, $l(m_x)$ is large)
is exotic in the sense that it belongs
to the kind of strings about which the halting sequence contains
much information and {\em vice versa}: $I(\chi : x)$ is large.
\end{remark}

 \begin{proof}
When we talk about complexity with $\chi$ in the condition, we use a 
Turing machine with $\chi$ as an ``oracle''.
With the help of $\chi$, we can compute $m_{x}$, and so we can define the
following new semicomputable (relative to $\chi$) function with
$c = 6/\pi^{2}$:
 \[
   \nu(x \mid \chi) = c \m(x) 2^{l(m_{x})} / l(m_{x})^{2}.
 \]
% Gacs: rewrote the proof below.
We have, using~\ref{t.X(r)-bd} and defining $Y(r) = X(r) \setminus X(r+1)$
so that $l(m_x)=r$ for $x \in Y(r)$:
 \[
 \begin{split}
 \sum_{x \in Y(r)} \nu(x \mid \chi) &=
 c r^{-2} 2^{r} \sum_{x \in Y(r)} 2^{-K(x)}
\\ &\le c r^{-2} 2^{r} 2^{-r+2} \le 4 c r^{-2}.
 \end{split}
 \]
Summing over $r$ gives $\sum_{x} \nu(x \mid \chi) \le 4$.
The theorem that $\m(x) = 2^{-K(x)}$ is maximal within multiplicative
constant among semicomputable semimeasures is also true relative to oracles.
Since we have established that
$\nu(x \mid \chi )/4$ is a semicomputable semimeasure, therefore
$\m(x \mid \chi) \gem \nu(x \mid \chi )$, or equivalently,
 \[
   K(x \mid \chi) \lea 
  -\log \nu(x \mid \chi ) \eqa K(x) - l( m_{x}) + 2 \log l(m_{x}),
 \]
which proves the theorem.
 \end{proof}

 \section{Non-Stochastic Objects}\label{s.non-stoch}

In this section, whenever we talk about a description of
a finite set $S$ we mean an {\em explicit} description. This
establishes the precise meaning of $K(S)$, $K(\cdot \mid S)$,
${\bf m} (S) = 2^{-K(S)}$, and ${\bf m} (\cdot \mid S) = 2^{K(\cdot \mid S)}$,
and so forth.

Every data sample consisting
of a finite string $x$ has an sufficient statistic in the form of
the singleton set $\{x\}$. Such a sufficient statistic is
not very enlightening since it simply replicates the data and
has equal complexity with $x$. 
Thus, one is interested in the minimal sufficient statistic
that represents the regularity, (the meaningful) information,
in the data and leaves out the accidental features. This raises
the question whether every $x$ has a minimal sufficient
statistic that is significantly less complex than
$x$ itself. At a Tallinn conference in 1973 Kolmogorov
(according to \cite{ShenStoch83,Co85}) raised
the question whether there are objects $x$ that
have no minimal sufficient statistic that have relatively 
small complexity. In other words, he inquired into the existence of
objects that are not in general position (random with respect to)
any finite set of small enough complexity, that is, 
``absolutely non-random'' objects. Clearly, such objects $x$ have
neither minimal nor maximal complexity: if they have minimal complexity then
the singleton set $\{x\}$ is a minimal sufficient statistic of small 
complexity, and
if $x \in \{0,1\}^n$
is completely incompressible (that is, it is individually random 
and has no meaningful information), then 
the uninformative universe $\{0,1\}^n$ is 
the minimal sufficient statistic of small complexity. To analyze
the question better we need the technical notion of randomness deficiency.

Define the {\em randomness deficiency} of an object $x$ with respect
to a finite set $S$ containing it as the amount by
which the complexity of $x$ as an element
of $S$ falls short of the maximal possible complexity
of an element in $S$ when $S$ is known explicitly (say, as a list):
 \begin{equation}\label{eq.randdef}
   \d_{S}(x) = \log | S| - K(x \mid S).
 \end{equation}
The meaning of this function is clear: most elements of $S$ have
complexity near $\log  |S|$, so this difference measures the amount of
compressibility in $x$ compared to the generic, typical, random elements of $S$.
This is a generalization of the sufficiency notion in that it
measures the discrepancy with typicality and hence sufficiency: 
if a set $S$ is a sufficient
statistic for $x$ then $\d_{S}(x) \eqa 0$.

We now continue the discussion of Kolmogorov's question.
Shen~\cite{ShenStoch83} gave a first answer by establishing the
existence of absolutely non-random objects $x$ of length $n$,
having randomness deficiency at least $n-2k - O(\log k)$ with respect
to every finite set $S$ of complexity $K(S) <k$ that
contains $x$. Moreover, since the set $\{x\}$ has complexity $K(x)$
and the randomness deficiency of $x$ with respect to this singleton set
is $\eqa 0$, it follows by choice of $k=K(x)$ that the complexity 
$K(x)$ is at least $n/2 - O(\log n)$.  

Here we sharpen this
result:
We establish the existence of absolutely non-random objects $x$ of length $n$,
having randomness deficiency at least $n-k$ with respect
to every finite set $S$ of complexity $K(S \mid n) <k$ that
contains $x$. Clearly, this is best possible since $x$
has randomness deficiency of at least $n-K(S \mid n)$ with every finite set $S$
containing $x$, in particular,  with complexity $K(S \mid n)$  more than a fixed
constant below $n$ the randomness deficiency exceeds that fixed constant.
That is, every sufficient statistic for $x$ has complexity at least $n$.
But if we choose $S= \{x\}$ then $K(S \mid n) \eqa K(x \mid n) \lea n$, and, moreover,
the randomness deficiency of $x$ with respect to $S$
is $n-K(S \mid n) \eqa 0$. Together
this shows that the absolutely nonrandom objects $x$ length $n$
of which we established
the existence have complexity $K(x \mid n) \eqa n$, and moreover,
 they have significant
randomness deficiency with respect to every set $S$ containing them
that has complexity significantly below their own complexity $n$.

\subsection{Kolmogorov Structure Function}

We first consider the relation between the minimal unavoidable randomness
deficiency of $x$ with respect to a set $S$ containing it,
when the complexity of $S$ is upper bounded by $\alpha$.
These functional relations are known as {\em Kolmogorov structure functions}.
Kolmogorov proposed a variant of the
function
 \begin{equation}\label{eq.kolm1}
   \h_{x}(\ag) = \min_{S} \setof{\log | S| : x\in S,\; K(S) < \ag},
\end{equation}
where $S \subseteq \{0,1\}^*$ is a finite set containing $x$,
the contemplated model for $x$, and $\ag$ is a nonnegative
integer value bounding the complexity of the contemplated $S$'s.
He did not specify what is meant by $K(S)$ but it was noticed
immediately, as the paper~\cite{ShenStoch99} points out, that the behavior
of $\h_{x}(\ag)$ is rather trivial if $K(S)$ is taken to be the complexity
of a program that lists $S$ without necessarily halting.
Section~\ref{s.implic} elaborates this point.
So, the present section refers to explicit descriptions only.

It is easy to see that for every increment $d$ we have
 \[
   \h_{x}(\ag + d) \le |\h_{x}(\ag) - d + O(\log d)|,
 \]
provided the right-hand side is non-negative, and 0 otherwise.
Namely, once we have an optimal set $S_{\ag}$ we can subdivide it in any
standard way into $2^{d}$ parts and take as $S_{\ag + d}$ the part
containing $x$.
Also, $\h_{x}(\ag) = 0$ implies $\ag \gea K(x)$, and, since the choice
of $S=\{x\}$ generally implies only $\ag \lea K(x)$ is meaningful
we can conclude $\ag \eqa K(x)$.
Therefore it seems better advised to consider the function
  \[
   \h_{x}(\ag) + \ag - K(x) 
     = \min_{S} \setof{\log | S | - (K(x) - \ag) : K(S) < \ag}
 \]
rather than (\ref{eq.kolm1}).
For technical reasons related to the later analysis, 
we introduce the following variant of
randomness deficiency (\ref{eq.randdef}):
 \[
   \d^{*}_{S}(x) = \log |S| - K(x \mid S, K(S)).
 \]
The function $\h_{x}(\ag) + \ag - K(x)$ seems related to a function
of more intuitive appeal, namely
$\bg_{x}(\ag)$ measuring the minimal unavoidable randomness deficiency
of $x$ with respect to every finite set $S$, that contains it,
of complexity $K(S) < \ag$.
Formally, we define
\[
   \bg_{x}(\ag) = \min_{S} \setof{\d_{S}(x) : K(S) < \ag},
\]
and its variant 
\[
   \bg^{*}_{x}(\ag) = \min_{S} \setof{\d^{*}_{S}(x) : K(S) < \ag},
\]
defined in terms of $\d^{*}_{S}$.
Note that $\bg_{x}(K(x)) \eqa \bg^*_{x}(K(x)) \eqa 0$.
These $\beta$-functions are related to, but different from,
the $\beta$ in (\ref{def.stoch}). 

To compare $\h$ and $\bg$, let us confine ourselves to binary strings of
length $n$.
We will put $n$ into the condition of all complexities.
  
  \begin{lemma}\label{t.h-vs-beta}
 $\bg^{*}_{x}(\ag \mid n) \lea \h_{x}(\ag \mid n) + \ag - K(x \mid n)$.
  \end{lemma}

 \begin{proof}
Let $S \ni x$ be a set with $K(S \mid n) \le \ag$ and 
assume $\h_{x}(\ag \mid n) = \log  |S|$.
Tacitly understanding $n$ in the conditions, and using the additivity
property (\ref{eq.soi}),
 \begin{align*}
    K(x) - \ag  \le K(x) - K(S)  & \lea K(x, S) - K(S) 
 \\&             \eqa K(x \mid S, K(S)).
 \end{align*}
Therefore
 \begin{align*}
   \h_{x}(\ag) + \ag - K(x) &= \log | S| - (K(x) - \ag) 
\\                          &\gea \log | S| - K(x \mid S, K(S))
                          \ge \bg^{*}_{x}(\ag). 
 \end{align*}
 \end{proof}

It would be nice to have an inequality also in the other direction, but we
do not know currently what is the best that can be said.

\subsection{Sharp Bound on Non-Stochastic Objects}
We are now able to formally express the notion of non-stochastic
objects using the Kolmogorov structure functions 
$\bg_{x}(\ag), \bg^*_{x}(\ag)$.
For every given $k < n$, Shen constructed in~\cite{ShenStoch83} a binary
string $x$ of length $n$ with $K(x) \le k$ and
$\bg_{x}(k - O(1)) > n - 2 k - O(\log k)$. 
Let $x$ be one of the non-stochastic objects of which
the existence is established. 
Substituting $k \eqa K(x)$  we can contemplate
the set $S= \{x\}$ with complexity $K(S) \eqa k$
and $x$ has randomness deficiency $\eqa 0$ with respect to $S$.
This yields $0 \eqa \bg_{x}(K(x)) \gea n - 2K(x) - O(\log K(x))$. 
Since it generally holds that these non-stochastic
objects have complexity $K(x) \gea n/2 - O(\log n)$,
they are {\em not random, typical, or in general position}
with respect to every set $S$ containing them with complexity
$K(S) \not\gea n/2 - O(\log n)$, but they are random, typical,
or in general position only for sets $S$ with complexity
$K(S)$ sufficiently exceeding  $n/2- O(\log n)$ like $S= \{x\}$.

% Clearly, we have $\eqa 0$
%minimal randomness deficiency for $k=K(x)$. Therefore,
%$0 \eqa \bg_{x}(K(x)) > n - 2 K(x) - O(\log K(x))$ which yields
%$K(x) \gea n/2 - O(\log n)$ (since $K(x) \lea n + 2 \log n$).

Here, we improve on this result, replacing $n - 2k - O(\log k)$ 
with $n - k$ and
using $\bg^*$ to avoid logarithmic terms.
This is the best possible, since by choosing $S = \{0,1\}^{n}$ we find
$\log  |S| - K(x \mid S, K(S)) \eqa n - k$, and hence
$\bg^{*}_{x}(c) \lea n - k$ for some constant $c$, which implies
$\bg^{*}_{x}(\ag) \le \bg_{x}(c) \lea n - k$ for every $\ag > c$.

 \begin{theorem}\label{t.non-stoch}
There are constants $c_{1}, c_{2}$ such that for any given $k < n$
there is a a binary
string $x$ of length $n$ with $K(x \mid n) \le k$ such that for all
 $\ag < k - c_{1}$ we have
 \[
  \bg^{*}_{x}(\ag \mid n) > n - k - c_{2}.
 \]
 \end{theorem}

In the terminology of (\ref{def.stoch}), 
the theorem states that there are constants $c_1,c_2$
such that for every $k < n$ there exists
a string $x$ of length $n$ of complexity $K(x \mid n ) \leq k$
that is not $(k-c_1, n-k-c_2)$-stochastic.

 \begin{proof}
Denote the conditional universal probability
as  $\m(S \mid n) = 2^{-K(S \mid n)}$.
We write ``$S \ni x$'' to indicate sets $S$ that satisfy $x \in S$.
For every $n$, let us define a function 
over all strings $x$ of length $n$ as follows:
 \begin{equation}\label{e.nu-def}
   \nu^{\le i}(x \mid n) = \sum_{S \ni x,\ K(S \mid n) \le i}
     \frac{\m(S \mid n)}{|S|}
 \end{equation}
% Gacs: new:
The following lemma shows that this function of $x$ is a semimeasure.

 \begin{lemma}
We have
 \begin{equation}\label{e.nu-meas}
   \sum_{x} \nu^{\le i}(x \mid n) \le 1.
 \end{equation}
 \end{lemma}
 \begin{proof}
We have 
 \[
  \begin{split}
  \sum_{x} \nu^{\le i}(x \mid n) 
    &\le \sum_{x}\sum_{S \ni x} \frac{\m(S \mid n)}{|S|}
     =  \sum_{S} \sum_{x \in S} \frac{\m(S \mid n)}{|S|}
\\  &=  \sum_{S} \m(S \mid n) \le 1.
  \end{split}
 \]
 \end{proof}
% Gacs: deleted the following:
% This $ \nu^{\le i}(x \mid n)$ approximates 
% $\m(x \mid n) = 2^{-K(x \mid n)}$ from below with growing $i$, since
% $K(\{x\} \mid n) \eqa K(x \mid n)$.

 \begin{lemma} There are constants $c_{1}, c_{2}$ such that for some $x$ of
length $n$,
  \begin{align}\label{e.x-prop}  
% Gacs: replaced \lem with \le:
        \nu^{\le k - c_{1}}(x \mid n ) &\le 2^{-n},
\\\label{e.K-x}
                             k - c_{2} &\le K(x \mid n) \le k.
  \end{align}
 \end{lemma}
 \begin{proof}
Let us fix $0 < c_{1} < k$ somehow, to be chosen appropriately later.
Inequality~(\ref{e.nu-meas}) implies that there is an $x$ 
with~(\ref{e.x-prop}).
Let $x$ be the first string of length $n$ with this property.
To prove the right inequality of (\ref{e.K-x}), let $p$ be the program
of length $\le i = k - c_{1}$
that terminates last in the standard 
running of all these programs simultaneously in dovetailed fashion,
on input $n$.
We can use $p$ and its length $l(p)$ to compute all programs of length
$\le l(p)$ that output finite sets using $n$.
This way we obtain a list of all sets $S$ with $K(S \mid n) \le i$.
Using this list, for each $y$ of length $n$ we can compute 
$\nu^{\le i}(y \mid n)$, by 
using the definition~(\ref{e.nu-def}) explicitly.
Since $x$ is defined as the first $y$ with 
$\nu^{\le i}(y \mid n) \le 2^{-n}$, we can thus find $x$
by using $p$ and some program of constant length.
If $c_{1}$ is chosen large enough, then this implies $K(x \mid n) \le k$.

On the other hand, from the definition~(\ref{e.nu-def}) we have
 \[
   \nu^{\le K(\{x\} \mid n)}(x \mid n) \ge 2^{- K(\{x\} \mid n)}.
 \]
This implies, by the definition of $x$, that either
 $K(\{x\} \mid n) > k - c_{1}$ or $K(\{x\} \mid n) \ge n$.
Since $K(x \mid n) \eqa K(\{x\} \mid n))$ we get the left inequality of
(\ref{e.K-x}) in both cases for an appropriate $c_{2}$.
 \end{proof}

Consider now a new semicomputable function
 \[
   \mu_{x, i}(S \mid n) = \frac{2^{n} \m(S \mid n)}{ |S|}
 \]
on all finite sets $S \ni x$ with $K(S \mid n) \le i$.
Then we have, with $i = k - c_{1}$:
% Gacs: supplied the requested proofs.
 \[
   \begin{split}
    \sum_{S} \mu_{x, i}(S \mid n) &= 
   2^{n} \sum_{S \ni x,\ K(S \mid n) \le i} \frac{\m(S \mid n)}{|S|}
\\ &= 2^{n} \nu^{\le i}(x \mid n) \le 1
   \end{split}
 \]
by~(\ref{e.nu-def}),~(\ref{e.x-prop}), respectively, and so
$\mu_{x, i}(S \mid n)$ with $x,i,n$ fixed is a
lower semicomputable semimeasure.
By the dominating property we have 
${\bf m}(S \mid x,i,n) \gem \mu_{x, i}(S \mid n)$. 
Since $n$ is the length of $x$ and $i \eqa k$ we can set
$K(S \mid x,i,n) \eqa K(S \mid x,k)$, and hence
$K(S \mid x, k)  \lea -\log \mu_{x,i}(S \mid n)$.
Then, with the first $\eqa$ because of ~(\ref{e.K-x}),
 \begin{equation}\label{e.S-on-x}
 \begin{split}
   K(S &\mid x, K(x \mid n))  
\\     &\eqa K(S \mid x, k)  \lea -\log \mu_{x,i}(S \mid n) 
\\     &= \log |S| - n +  K(S \mid n).
 \end{split}
 \end{equation}
Then, by the additivity property~(\ref{eq.soi}) and~(\ref{e.S-on-x}): 
 \begin{align*}
     K(x &\mid S, K(S \mid n), n)
\\       &\eqa K(x \mid n) + K(S \mid x, K(x \mid n)) - K(S \mid n) 
\\       &\lea k + \log |S| - n.
 \end{align*}
Hence 
$\d^{*}(x \mid S, n) = \log |S| - K(x \mid S, K(S \mid n), n) \gea n - k$.
 \end{proof}

We are now in the position to prove Lemma~\ref{lem.emss}:
For every length $n$, there exist strings $x$ of length $n$
with $K(x \mid n) \eqa n$ for which $\{x\}$ is an explicit
minimal sufficient statistic.

\begin{proof} ({\em of Lemma~\ref{lem.emss}}):
Let $x$ be one of the non-stochastic objects of which
the existence is established by Theorem~\ref{t.non-stoch}.
Choose $x$ with $K(x \mid n) \eqa k$  so that the set
$S= \{x\}$ has complexity $K(S \mid n) = k-c_1$
and $x$ has randomness deficiency $\eqa 0$ with respect to $S$. 
Because $x$ is non-stochastic,
this yields $0 \eqa \bg^{*}_{x}(k-c_1 \mid n) \gea n-K(x \mid n)$. 
For every $x$ we have $K(x \mid n) \lea n$. Together it follows that
$K(x \mid n) \eqa n$. That is, these non-stochastic
objects $x$ have complexity $K(x \mid n) \eqa n$. Nonetheless,
there is a constant $c'$ such that $x$
is {\em not random, typical, or in general position}
with respect to any explicitly represented finite set $S$ 
containing it that has complexity
$K(S \mid n)  < n-c'$, but they are random, typical, 
or in general position for some sets $S$ with complexity
$K(S \mid n) \gea n$ like $S= \{x\}$.
That is, every explicit sufficient statistic $S$ for $x$ 
has complexity $K(S \mid n) \eqa n$, and $\{x\}$ is
such a statistic. Hence $\{x\}$ is an explicit minimal sufficient
statistic for $x$.
\end{proof}

\section{Probabilistic Models}
\label{s.prob}
It remains to generalize the model class from finite sets to 
the more natural and significant setting of probability distributions.
Instead of finite sets the models are computable probability
density functions $P: \{0,1\}^* \rightarrow [0,1]$ with
$\sum P(x) \leq 1$---we allow defective probability distributions
where we may concentrate the surplus probability on a 
distinguished ``undefined'' element.
``Computable'' means that there is a Turing machine $T_P$ that computes
approximations to 
the value of $P$ for every argument
(more precise definition follows below).
The (prefix-) complexity $K(P)$ of a 
computable partial function $P$ is defined by 
\[
K(P) = \min_i \{K(i): \mbox{\rm Turing machine } T_i 
\; \; \mbox{\rm computes }
P \}. 
\]
Equality~(\ref{eq.typ}) now becomes
\begin{equation}\label{eq.typP}
K(x \mid P^*) \eqa - \log P(x) ,
\end{equation}
and equality~(\ref{eq.optim}) becomes
\[
K(x) \eqa K(P)  - \log P(x) .
\]
As in the finite set case, the complexities involved are crucially
dependent on what we mean by ``computation'' of $P(x)$, that is,
on the requirements on the format in which the output is to be represented.
Recall from \cite{LiVi97} that Turing machines can compute rational numbers:
If a Turing machine $T$ computes $T(x)$, then we interpret
the output as a pair of natural numbers,
$T(x) = \langle p,q \rangle$, according to a standard pairing function. 
Then, the rational value computed by $T$ is by definition $p/q$.
The distinction between explicit and implicit description of $P$
corresponding to the finite set model case is now defined as follows:
\begin{itemize}
\item
It is {\em implicit} if  there
is a Turing machine $T$ computing $P$
halting with rational value $T(x)$ so that
$- \log T(x) \eqa - \log P(x)$, and,
furthermore, $K(- \log T(x) \mid P^*) \eqa 0$ for $x$
satisfying (\ref{eq.typP})---that is, for typical $x$.
\item
It is {\em explicit} if the Turing machine $T$ computing $P$, given
$x$ and a tolerance $\epsilon$ halts with rational value
so that
$- \log T(x) = - \log ( P(x) \pm \epsilon)$, and,
furthermore, $K(- \log T(x) \mid P^*) \eqa 0$ for $x$
satisfying (\ref{eq.typP})---that is, for typical $x$.
%\item
%It is {\em implicit} if there are positive
%constants $c,C$ such that the Turing machine $T$ computing $P$
%halts with rational value $T(x)$ with $cP(x) < T(x) < C P(x)$.
%Hence $- \log T(x) \eqa - \log P(x)$.
%\item
%It is {\em explicit} if the Turing machine $T$ computing $P$, given
%$x$ and a tolerance $\epsilon$ halts with rational value 
%$P(x) - \epsilon < T(x) < P(x) + \epsilon$.
\end{itemize}
The implicit and explicit descriptions of finite sets and
of uniform distributions with $P(x) = 1/|S|$ for all $x \in S$
and $P(x)=0$ otherwise, are as follows:
An implicit (explicit) description
of $P$ is identical with an implicit (explicit) description
of $S$, up to a short fixed program which indicates which of the two
is intended, so that $K(P(x)) \eqa K(S)$ for $P(x) > 0$ (equivalently,
$x \in S$).

To complete our discussion: the worst case of representation format,
a recursively enumerable approximation of $P(x)$
where nothing is known about its value,
would lead to indices $- \log P(x)$ of unknown length.
We do not consider this case.

The properties for the probabilistic
models are loosely related to the properties of 
finite set models by Proposition~\ref{prop.1}.
We sharpen the relations by appropriately modifying
the treatment of the finite set case, but essentially following the same
course.
 
We may use the notation
 \[
   P_{\impl}, P_{\expl}
 \]
for some implicit and some explicit representation of $P$.
When a result applies to both implicit and explicit representations,
or when it is clear from the context which representation is meant, we
will omit the subscript.

\subsection{Optimal Model and Sufficient Statistic}
As before, we distinguish between ``models'' that are
computable probability distributions,
and the ``shortest programs'' to compute those models
that are finite strings.

Consider a string $x$
of length $n$ and prefix complexity $K(x)=k$.
We identify the {\em structure} or {\em regularity} in $x$ that are
to be summarized with a computable probability density function $P$
with respect to which $x$ is a {\em random} or  {\em typical} member.
%The shortest binary program $P^*$ for $P$ is called an
%{\em algorithmic statistic} for data $x$.
For $x$ typical for $P$ holds the following \cite{LiVi97}:
Given 
an (implicitly or explicitly described)
shortest program $P^*$ for $P$, a shortest binary program
computing $x$ (that is, of length $K(x \mid P^*)$)
can not
be significantly shorter
than its Shannon-Fano code \cite{CT91} of
length $- \log P(x)$, that is,
$ K(x \mid P^*) \gea - \log P(x) $.
%pg:
By definition, we fix some agreed upon constant
 $
  \beta \ge 0,
 $
and require 
\[ K(x \mid P^*) \ge - \log P(x)  - \beta.
\]
As before, we will not indicate the dependence on $\beta$ explicitly, but the
constants in all our inequalities ($\lea$) will be allowed to be functions
of this $\beta$.
This definition requires a positive $P(x)$.
In fact, since
$K(x \mid P^*) \lea K(x) $, it limits the size of $P(x)$ to 
$\Omega (2^{-k})$.
The shortest program $P^*$ from which a probability density 
function $P$ 
can be computed is an {\em algorithmic statistic} for $x$ iff
\begin{equation}
 K(x \mid P^*) \eqa - \log P(x) .
\end{equation}
There are two natural measures of suitability of such a statistic.
We might prefer either the simplest 
distribution, or the largest distribution, as
corresponding to the most likely structure `explaining' $x$.
The singleton probability distribution $P(x)=1$,
while certainly a statistic for $x$,
would indeed be considered a poor explanation.
Both measures relate to the optimality of a two-stage description of
$x$ using $P$:
\begin{align}\label{eq.twostageP}
 K(x) \leq K(x,P) & \eqa  K(P) + K(x \mid P^*) 
\\& \lea K(P) - \log P(x),
\nonumber
\end{align}
where we rewrite $K(x,P)$ by~(\ref{eq.soi}).
Here, $P$ can be understood as either $P_{\impl}$ or $P_{\expl}$.
Call a distribution $P$ (with positive probability $P(x)$)
 for which
\begin{equation}\label{eq.optimP}
K(x) \eqa K(P) - \log P(x),
\end{equation}
{\em optimal}.
(More precisely, we should require $K(x) \ge K(P) - \log P(x) - \beta$.)
Depending on whether $K(P)$ is understood as $K(P_{\impl})$ or
$K(P_{\expl})$, our
definition splits into implicit and explicit optimality.
The shortest program for an optimal computable probability distribution
is a {\em algorithmic sufficient statistic} for $x$.

\subsection{Properties of Sufficient Statistic}

As in the case
of finite set models , we start with a sequence of lemmas 
that are used to obtain the main results on minimal sufficient
statistic.
Several of these lemmas have two versions: for implicit distributions
 and for explicit
distributions.
In these cases, $P$ will denote $P_{\impl}$ or $P_{\expl}$ respectively.

Below it is shown that the mutual information between every
typical distribution and the
data is not much less than $K(K(x))$, the complexity of the complexity
$K(x)$ of the data $x$.
For optimal distributions it is at least that, and for
algorithmic minimal statistic it is equal to that.
%pg: I do not understand the following sentence.
% This suggests that the mutual information must be larger
% for more ``random'' data complexities.
The log-probability
of a typical distribution is determined by the following:

\begin{lemma}\label{lem.cardP}
Let $k=K(x)$.
If a distribution $P$ is
(implicitly or explicitly)
typical for $x$ then
$I(x:P) \eqa k + \log P(x)$.
%$|S| = \Theta (2^{k-I(x:S)})$.
\end{lemma}
\begin{proof}
By definition $I(x:P) \eqa K(x) - K(x \mid P^*)$
and by typicality $K(x \mid P^*) \eqa - \log P(x)$.
\end{proof}

The above  lemma states that for
(implicitly or explicitly)
typical $P$ the probability $P(x) = \Theta ( 2^{- (k-I(x:P))})$.
The next lemma asserts that for implicitly typical $P$
the value $I(x:P)$ can fall below $K(k)$ by no more than an
additive logarithmic term.
% plus the amount of information
%required to compute $- \log P(x)$ from $P$.

\begin{lemma}\label{lem.KIP}
Let $k=K(x)$.
If a distribution $P$ is (implicitly or explicitly)
typical for $x$ then $I(x:P) \gea K(k) - K(I(x:P))$
%typical for $x$ then $I(x:P) \gea K(k) - K(I(x:P)) - K(- \log P(x) \mid P^*)$
and $- \log P(x) \lea k-  K(k) + K(I(x:P))$.
%and $- \log P(x) \lea k-  K(k) + K(I(x:P)) + K(- \log P(x) \mid P^*)$.
(Here, $P$ is understood as $P_{\impl}$ or $P_{\expl}$ respectively.)
\end{lemma}

\begin{proof}
Writing $k=K(x)$, since
\begin{equation}\label{eq.kxP}
k \eqa K(k,x) \eqa K(k)+K(x \mid k^*)
\end{equation}
 by~(\ref{eq.soi}), we have
$I(x:P) \eqa K(x)-K(x \mid P^*) \eqa K(k)-[K(x \mid P^*)-K(x \mid k^*)]$.
Hence, it suffices to show $K(x \mid P^*)-K(x \mid k^*) \lea
K(I(x:P))$. %  + K(- \log P(x) \mid P^*)$.
Now, from an implicit description $P^*$
% and a program $q$
%of length $\eqa  K(- \log P(x) \mid P^*)$ we can
we can 
find the value $\eqa - \log P(x) \eqa k-I(x:P)$. To recover $k$ from $P^*$,
we at most
require an extra $K(I(x:P))$ bits.
That is, $K(k \mid P^*) \lea K(I(x:P))$. % +  K(- \log P(x) \mid P^*)$.
This reduces what we have to show to
$K(x \mid P^*) \lea K(x \mid k^*) + K(k \mid P^*)$ which is asserted by
Theorem~\ref{lem.magic}. This shows the first statement in the theorem.
The second statement follows from the first one: rewrite
$I(x:P) \eqa k + K(x \mid P^*)$
and substitute $- \log P(x) \eqa K(x \mid P^*)$.
\end{proof}

If we further restrict typical distributions to optimal ones then
the possible positive probabilities assumed
by distribution $P$ are slightly restricted.
First  we show that implicit optimality with respect to 
some data
 is equivalent to typicality with respect to the data
combined with effective constructability (determination) from the data.

\begin{lemma}\label{lem.iotyP}
A distribution $P$ is
(implicitly or explicitly)
optimal for $x$ iff it is typical and $K(P \mid x^*) \eqa 0$.
\end{lemma}

\begin{proof}
A distribution $P$ is optimal iff~(\ref{eq.twostageP}) holds with equalities.
Rewriting $K(x,P) \eqa K(x)+K(P \mid x^*)$
the first inequality becomes an equality iff $K(P \mid x^*) \eqa 0$,
and the second inequality becomes an equality iff $K(x \mid P^*)
\eqa - \log P(x)$
(that is, $P$ is a typical distribution).
\end{proof}

\begin{lemma}\label{lem.optP}
Let $k=K(x)$.
If a distribution $P$ is (implicitly or explicitly) optimal for $x$, then
$I(x:P) \eqa K(P) \gea K(k)$. %- K(- \log P(x) \mid P^*)$
and $- \log P(x) \lea  k-K(k)$. %+K(- \log P(x) \mid P^*)$.
\end{lemma}

\begin{proof}
If $P$ is optimal for $x$, then
$k = K(x) \eqa K(P)+K(x \mid P^*) \eqa K(P)-\log P(x)$.
%It follows that $K(k \mid S^*) \eqa K(K(S)+(\log |S|) \mid S^*) \eqa 0$,
%hence $K(k) \lea K(S) + K(k \mid S^*) \eqa K(S)$.
 From $P^*$ %and a program $q$ of length $K(- \log P(x) \mid P^*)$,
we can find both $K(P) \eqa l(P^*)$
and $ \eqa - \log P(x)$, and hence $k$, that is,
$K(k) \lea K(P)$. %+K(- \log P(x)\mid P^*)$.
We have $I(x:P) \eqa K(P)-K(P \mid x^*) \eqa K(P)$ by~(\ref{eq.soi}),
Lemma~\ref{lem.iotyP}, respectively.
This proves the first property.
Substitution of $I(x:P) \gea K(k)$
% - K(- \log P(x) \mid P^*) $ 
in the expression
of 
Lemma~\ref{lem.cardP} proves the second property.
\end{proof}

\begin{remark}
Our definitions of implicit and explicit description format
entail that, for typical $x$,
one can compute $\eqa - \log P(x)$ and $- \log P(x)$,
respectively, from $P^*$ alone without requiring $x$. An alternative
possibility would have been that implicit and explicit description
formats refer to the fact that we can compute $\eqa - \log P(x)$
and $- \log P(x)$, respectively, given {\em both} $P$ and $x$.
This would have added a $-K(- \log P(x) \mid P^*)$ additive
term in the righthand side of the expressions 
in Lemma~\ref{lem.KIP} and Lemma~\ref{lem.optP}.
Clearly, this alternative definition is equal to
the one we have chosen  iff this term is always $\eqa 0$
for typical $x$. We now show that this is not the case. 

Note that for distributions that are uniform
(or almost uniform) on a finite support
we have $K(- \log P(x) \mid P^*) \eqa 0$: In this borderline
case the result specializes to that of Lemma~\ref{lem.KI} for finite
set models, and the two possible definition types for implicitness 
and those for explicitness coincide.

On the other end of the spectrum, for the definition type considered
in this remark,
the given lower bound on $I(x:P)$ drops
in case knowledge of $P^*$
doesn't suffice to compute $- \log P(x)$, that is, if
$K ( - \log P(x) \mid P^*) \gg 0$ for an statistic $P^*$ for $x$.
The question is, whether we can exhibit
such a probability distribution that is also computable?
The answer turns out to be affirmative.
By a result due to R. Solovay and P. G\'acs, \cite{LiVi97} Exercise
3.7.1 on p. 225-226, there is a computable function $f(x) \gea K(x)$
such that
$f(x) \eqa K(x)$ for infinitely many $x$.
Considering the case of $P$
optimal
for $x$ (a stronger assumption than that $P$ is just typical) we have
$- \log P(x) \eqa K(x) - K(P)$.
Choosing
$P(x)$ such that $- \log P(x) \eqa  \log f(x) - K(P)$, we have that
$P(x)$ is computable since $f(x)$ is computable and $K(P)$ is
a fixed constant. Moreover, there
are infinitely many $x$'s for which $P$ is optimal, so
$K(- \log P(x) \mid P^*) \rightarrow \infty$ for $x \rightarrow \infty$
through this special sequence.
\end{remark}

%Note that for distributions that are uniform
%(or almost uniform) on a finite support
%we have $K(- \log P(x) \mid P^*) \eqa 0$: In this borderline
%case the result specializes to that of Lemma~\ref{lem.opt} for finite
%set models.

%On the other end of the spectrum,
%we have the case that knowledge of $P^*$
%doesn't help to compute $- \log P(x)$, that is,
%$K(- \log P(x) \mid P^*) \gg 0$ as exemplified
%above.
%Then, the lower bound on $I(x:P) \eqa K(P)$ drops 
%towards 0 while
%the upper bound on $- \log P(x)$ rises towards $k$.

\subsection{Concrete Minimal Sufficient Statistic}\label{s.implicP}

A simplest implicitly optimal distribution (that is, of least complexity)
is an implicit algorithmic minimal sufficient statistic. As before,
let $S^k = \{y: K(y) \leq k \}$.  
Define the distribution $P^k (x) = 1/|S^k|$ for $x \in S^k$, and
$P^k (x)= 0$ otherwise.
The demonstration that $P^k(x)$
is an implicit algorithmic minimal sufficient statistic proceeeds
completely analogous to the finite set model setting, 
Corollary~\ref{corr.iamss}, using the substitution
$K(- \log P^k(x) \mid (P^k)^*) \eqa 0$.

A similar equivalent construction suffices to obtain an explicit
algorithmic minimal near-sufficient statistic for $x$,
analogous to $S^k_{m_x}$ in the finite set model setting,
Theorem~\ref{theo.ekmss}.
That is, $P^k_{m_x} (y) = 1/|S^k_{m_x}|$ for $y \in S^k_{m_x}$,
and 0 otherwise. 

In general, one can develop the theory of minimal sufficient statistic
for models that are probability distributions similarly to that
of finite set models.
%, up to the extra additive term
%$K(- \log P(x) \mid P^*)$. It is not known how far that term
%can be reduced.

\subsection{Non-Quasistochastic Objects}
As in the more restricted case of finite sets, there are objects
that are not typical for any explicitly computable probability distribution
that has complexity significantly below that of the object itself.
With the terminology of (\ref{def.quasistoch}), we may call such
{\em absolutely non-quasistochastic}.
%First, adapt the relevant definitions $\d^{*}_{S}(x)$ and
%$\bg_{x}^{*} (\ag)$ related
%to explicitly represented finite sets to
%explicitly computable probability distributions as $\d^{*}_{P}(x)$
%and $\pi_{x}^{*} (\ag)$, respectively:
 %\begin{eqnarray*}
   %\d^{*}_{P}(x) & = & - \log P(x) - K(x \mid P, K(P)) \\
 %\pi_{x}^{*} (\ag) & = & \min_{P} \setof{\d^{*}_{P}(x) : K(P) < \ag}
 %\end{eqnarray*}

By Proposition~\ref{prop.1}, item (b), there are constants $c$ and $C$
such that if $x$ is
not $(\alpha + c \log n, \beta + C)$-stochastic (\ref{def.stoch})
then $x$ is not $(\alpha, \beta )$-quasistochastic (\ref{def.quasistoch}). 
Substitution in Theorem~\ref{t.non-stoch} yields:

 \begin{corollary}\label{t.non-stochP}
There are constants $c,C$ such that,
for every $k < n$, there are constants $c_{1}, c_{2}$ and a binary
string $x$ of length $n$ with $K(x \mid n) \le k$ such that 
$x$ is not $(k-c \log n - c_1 , n-k-C-c_2)$-quasistochastic.
 \end{corollary}

As a particular consequence: 
Let $x$ with length $n$  be one of the non-quasistochastic strings of which
the existence is established by Corollary~\ref{t.non-stochP}.
Substituting $ K(x\mid n) \lea k - c \log n $,
we can contemplate
the distribution $P_x(y)= 1$ for $y=x$ and 
and $0$ otherwise. Then we have
complexity $K(P_x\mid n) \eqa K(x \mid n)$.
Clearly, $x$ has randomness deficiency $\eqa 0$ with respect to $P_x$.
Because of the assumption of non-quasistochasticity of $x$,
and because the minimal randomness-deficiency $\eqa n-k$
of $x$ is always nonnegative, 
 $0 \eqa n - k  
\gea n-K(x \mid n) - c \log n$.
Since it generally holds that $K(x \mid n ) \lea n$, it follows that
$n \gea K(x \mid n) \gea n - c \log n$. That is, these non-quasistochastic
objects have complexity $K(x \mid n) \eqa n - O( \log n)$
and are {\em not random, typical, or in general position}
with respect to any explicitly computable distribution
 $P$ with $P(x) > 0$ and  complexity
$K(P \mid n) \lea n - (c+1) \log n$, but they are random, typical,
or in general position only for some distributions $P$ with complexity
$K(P \mid n) \gea n- c \log n$ like $P_x$.
That is, every explicit sufficient statistic $P$ for $x$
has complexity $K(P \mid n) \gea n - c \log n$, and $P_x$ is
such a statistic.

\section{Algorithmic Versus Probabilistic}
\label{sect.formanal}
Algorithmic sufficient statistic, a function of the data,
is so named because intuitively
it expresses an individual summarizing of the relevant information
in the individual data, reminiscent of 
the probabilistic sufficient statistic that summarizes the
relevant information in a data random variable about a model
random variable. Formally, however, previous authors have
not established any relation. Other algorithmic notions
have been successfully related to their probabilistic
counterparts. The most significant one is that for every computable
probability distribution, the expected prefix complexity of the
objects equals the entropy of the distribution up to an additive
constant term, related to the complexity of the distribution in
question. We have used this property in (\ref{eq.eqamipmi})
to establish a similar relation between the expected
algorithmic mutual information and the probabilistic mutual information.
We use this in turn to show that 
there is a close relation between the algorithmic version and
the probabilistic version of sufficient
statistic: A probabilistic sufficient statistic is 
with high probability a natural conditional form 
of algorithmic sufficient statistic  
for individual data, and, conversely, that with
high probability a natural conditional
form of algorithmic sufficient statistic is  also a probabilistic
sufficient statistic.

Recall the terminology of probabilistic mutual information
(\ref{eq.mutinfprob})
and probabilistic sufficient statistic (\ref{eq.suffstatprob}).
Consider a probabilistic ensemble of models,
a family of computable probability mass functions $\{f_{\theta} \}$
indexed by a discrete parameter $\theta$, together with a computable 
distribution $p_1$ over $\theta$.
(The finite set model case is the  restriction where
the $f_{\theta}$'s are restricted to uniform distributions
with finite supports.)
This way we have a random variable $\Theta$ with outcomes in $\{f_{\theta} \}$
and a random variable $X$ with outcomes 
in the union of domains of $f_{\theta}$, and 
$p(\theta,x) = p_1 (\theta) f_{\theta}(x)$ is computable.

\begin{notation}
To compare the algorithmic sufficient statistic
with the probabilistic sufficient statistic it is
convenient to denote the sufficient statistic
 as a function $S(\cdot)$ of the data in both cases. 
Let a statistic
$S(x)$ of data $x$ be the more general form of probability distribution
as in Section~\ref{s.prob}. That is, $S$ maps the data $x$ to the 
parameter $\rho$ that determines 
a probability mass function $f_{\rho}$ (possibly not an element
of $\{f_{\theta} \}$). Note that ``$f_{\rho} (\cdot)$'' corresponds
to ``$P(\cdot)$''
in Section~\ref{s.prob}.
If $f_{\rho}$ is computable, then this can be the 
Turing machine $T_{\rho}$ that computes
$f_{\rho}$.
Hence, in the current section, 
``$S(x)$'' denotes a probability distribution, say $f_{\rho}$,
and ``$f_{\rho}(x)$'' is the probability $f_{\rho}$ concentrates on data $x$.
\end{notation}
\begin{remark}
In the probabilistic statistics setting,
Every function $T(x)$ is a statistic of $x$, but only some
of them are a sufficient statistic. In the algorithmic statistic
setting we have a quite similar situation. In the finite set statistic
case $S(x)$ is a finite set, and in the computable probability
mass function case $S(x)$ is a computable probability mass function. 
In both algorithmic cases we have shown $K(S(x) \mid x^*) \eqa 0$
for $S(x)$ is an implicitly or explicitly described sufficient statistic. 
This means that the number of such sufficient statistics for $x$
is bounded by a universal constant, and that there is a universal program
to compute all of them from $x^*$---and hence to compute
the minimal sufficient statistic from $x^*$.
\end{remark}
\begin{lemma}\label{theo.eqpral}
Let $p(\theta,x) = p_1 (\theta) f_{\theta} (x)$ be a computable joint
probability mass function, and let
$S$ be a function. Then all three conditions below are equivalent
and imply each other:

(i) $S$ is a probabilistic sufficient statistic 
(in the form $I(\Theta, X) \eqa I(\Theta , S(X))$). 

(ii) $S$ satisfies
\begin{equation}\label{eq.eqami}
\sum_{\theta,x} p(\theta,x) I(\theta:x)
\eqa 
\sum_{\theta,x} p(\theta,x) I(\theta: S(x))
\end{equation}

(iii) $S$ satisfies
\begin{align*}
I(\Theta ; X) \eqa I(\Theta ; S(X)) & \eqa
\sum_{\theta,x} p(\theta,x) I(\theta:x) 
\\& \eqa
\sum_{\theta,x} p(\theta,x) I(\theta: S(x)).
\end{align*}

All $\eqa$ signs hold up to an $\eqa \pm 2K(p)$ constant additive term.

\end{lemma}

\begin{proof}
Clearly, (iii) implies (i) and (ii).

We show that both (i) implies (iii) and (ii) implies (iii):
By (\ref{eq.eqamipmi}) we have 
\begin{align}\label{eq.asseq}
I(\Theta ; X) & \eqa \sum_{\theta,x} p(\theta,x) I(\theta:x),
\\ I(\Theta ; S(X)) & \eqa \sum_{\theta,x} p(\theta,x) I(\theta: S(x)),
\nonumber
\end{align}
where  we absorb a $\pm 2K(p)$ additive term in the $\eqa$ sign.
Together with (\ref{eq.eqami}),
(\ref{eq.asseq}) implies 
\begin{equation}\label{eq.eqpmi}
 I(\Theta ; X) \eqa I(\Theta ; S(X)) ;
\end{equation}
and {\em vice versa} (\ref{eq.eqpmi}) together with (\ref{eq.asseq})
implies (\ref{eq.eqami}). 

\end{proof}

\begin{remark}
It may be worth stressing that $S$ in Theorem~\ref{theo.eqpral} can
be any function, without restriction.
\end{remark}

\begin{remark}
Note that (\ref{eq.eqpmi}) involves equality $\eqa$
rather than precise equality as in the
definition of the probabilistic sufficient
statistic (\ref{eq.suffstatprob}).
\end{remark}

\begin{definition}\label{def.thetaI}
Assume the terminology and notation above. 
A statistic $S$ for data $x$ 
is {\em $\theta$-sufficient with deficiency $\delta$} 
if 
$I(\theta , x) \eqa I(\theta , S(x)) + \delta$.
If $\delta \eqa 0$ then $S(x)$ is simply a {\em $\theta$-sufficient
statistic}.
\end{definition}

The following lemma shows that $\theta$-sufficiency is a type
of conditional sufficiency:

\begin{lemma}\label{claim.1}
Let $S(x)$ be a sufficient statistic for $x$. Then,
\begin{equation}\label{eq.theta}
 K(x \mid \theta^*) + \delta \eqa  K(S(x) \mid \theta^* ) - \log S(x).
\end{equation}
iff $I(\theta , x) \eqa I(\theta , S(x)) + \delta$.
\end{lemma}

\begin{proof}
(If) By assumption,
 $K(S(x)) - K(S(x) \mid  \theta^*) + \delta \eqa K(x) - K(x \mid \theta^*)$.
%that is, $I(\theta: S(x)) + \delta \eqa I(\theta:x)$. 
%Since $S$ is a sufficient statistic for $x$, the term
Rearrange and add 
$-K(x \mid S(x)^*)- \log S(x) \eqa 0$ (by typicality)
 to the right-hand side to obtain
$K(x \mid \theta^*) +K(S(x)) \eqa K(S(x) \mid \theta^*) + K(x)
- K(x \mid S(x)^*) - \log S(x) - \delta$.
Substitute according to $K(x) \eqa K(S(x))+K(x \mid S(x)^*)$
(by sufficiency) in the
right-hand side, and subsequently subtract 
$K(S(x))$ from both sides, to obtain
(\ref{eq.theta}).

(Only If) Reverse the proof of the (If) case. 

\end{proof}

The following theorems state that $S(X)$ is a probabilistic sufficient
statistic iff $S(x)$ is an algorithmic $\theta$-sufficient statistic,
up to small deficiency, with high probability.

\begin{theorem}
Let $p(\theta,x) = p_1 (\theta) f_{\theta} (x)$ be a computable joint
probability mass function, and let
$S$ be a function.
If $S$ is 
a recursive probabilistic sufficient statistic, then
$S$ is 
a $\theta$-sufficient statistic with deficiency $O(k)$,
with $p$-probability at least $1 - \frac{1}{k}$.
\end{theorem}

\begin{proof}
If $S$ is a probabilistic sufficient statistic,
then, by Lemma~\ref{theo.eqpral}, equality of $p$-expectations (\ref{eq.eqami})
holds. However, it is still consistent with this to have
large positive and negative differences
$I(\theta: x) -I(\theta:S(x))$
for different $(\theta,x)$ arguments, such that these
differences cancel each other. 
This problem is resolved by appeal to
the algorithmic mutual information non-increase
law (\ref{eq.nonincrease}) which shows that all differences are 
essentially positive:
$I(\theta : x) - I(\theta : S(x)) \gea -K(S)$.
Altogether, let $c_1,c_2$ be least positive constants such that
$I(\theta : x) - I(\theta : S(x))+c_1$ is always nonnegative
and its $p$-expectation is $c_2$.
Then, by Markov's inequality, 
\[
p ( I( \theta : x) - I(\theta : S(x)) \geq kc_2 - c_1 ) \leq \frac{1}{k},
\]
that is,
\[ p ( I( \theta : x) - I(\theta : S(x)) < kc_2 - c_1 ) 
> 1 - \frac{1}{k}.
\]
\end{proof}

\begin{theorem}
For each $n$, consider the set of data $x$ of length $n$.
Let $p(\theta,x) = p_1 (\theta) f_{\theta} (x)$ be a computable joint
probability mass function, and let
$S$ be a function.
If $S$ is an algorithmic $\theta$-sufficient statistic for 
$x$, with $p$-probability
at least $1-\epsilon$ ($1/\epsilon \eqa n + 2 \log n$), then
$S$ is a probabilistic sufficient statistic.
\end{theorem}

\begin{proof}
By assumption, using Definition~\ref{def.thetaI},
there is a positive constant $c_1$, such that, 
\[
p ( | I(\theta : x) - I(\theta : S(x))| \leq c_1) \geq 1- \epsilon.
\]
Therefore,
\begin{align*}
0 \leq \sum_{| I(\theta : x )  - I(\theta : S(x))| \leq  c_1 } p(\theta ,x) 
& |I(\theta : x )  - I(\theta : S(x))|
\\ &  \lea  (1-\epsilon)c_1 \eqa  0. 
\end{align*}
On the other hand, since
 \[
1/\epsilon \gea n + 2 \log n \gea K(x) \gea  \max_{\theta , x} I(\theta ; x),
\]
we obtain
\begin{align*}
0 \leq \sum_{| I(\theta : x )  - I(\theta : S(x))| >  c_1 } p(\theta ,x) 
& |I(\theta : x )  - I(\theta : S(x))|
\\ &  \lea  \epsilon (n+2 \log n) \lea  0. 
\end{align*}
Altogether, this implies (\ref{eq.eqami}), and by
Lemma~\ref{theo.eqpral}, the theorem.
\end{proof}

\section{Conclusion}
\label{s.discussion}
An algorithmic
sufficient statistic is 
an individual finite set (or probability distribution) for which a given
individual sequence is a typical member. The theory is formulated
in Kolmogorov's absolute notion of the quantity of information in an 
individual object. 
This is a notion analogous to, and in some sense sharper than the
probabilistic  notion
of sufficient statistic---an average notion based on
the entropies of random variables. It turned out, that for every
sequence $x$ we can determine the complexity range of possible
algorithmic sufficient statistics,
and, in particular, exhibit a algorithmic minimal sufficient statistic.
The manner in which the statistic is effectively represented
is crucial: we distinguish implicit representation and
explicit representation. The latter is essentially a list
of the elements of a finite set or a table of the probability 
density function; the former is less explicit than a list
or table but more explicit than just recursive enumeration
or approximation in the limit. The algorithmic minimal sufficient statistic
can be considerably more complex depending on whether we want
explicit or implicit representations. We have shown that there are sequences
that have no simple explicit algorithmic sufficient statistic: the
algorithmic  minimal
sufficient statistic is essentially the sequence itself. Note that
such sequences cannot be random in the sense of having maximal
Kolmogorov complexity---in that case already the simple set of
all sequences of its length,
or the corresponding uniform distribution,
is an algorithmic  sufficient statistic of almost zero complexity.
We demonstrated close relations between the probabilistic
notions and the corresponding algorithmic notions:
(i) The average algorithmic mutual information
is equal to the probabilistic mutual information.
(ii)  
To compare algorithmic sufficient statistic
and probabilistic sufficient statistic meaningfully
one needs to consider a conditional version of algorithmic sufficient
statistic. We defined such a notion and demonstrated that
probabilistic sufficient statistic is with high probability an
(appropriately conditioned) algorithmic sufficient statistic and vice versa.
The most conspicuous theoretical open end is as follows:
For explicit descriptions we were only able to guarantee
a algorithmic minimal near-sufficient statistic, although the construction
can be shown to be minimal sufficient for almost all sequences.
One would like to obtain a concrete example of a truly 
explicit algorithmic minimal sufficient statistic.
%In the theory of sufficient statistic
%for models that are probability distributions, in contrast to that
%of finite set models, one has to deal with an extra additive term
%$K(- \log  P(x) \mid P^*)$. It is not known how far that term
%can be reduced.

\subsection{Subsequent Work}
One can continue generalization of model classes for algorithmic
statistic beyond computable probability mass functions. The ultimate
model class is the set of recursive functions. In the manuscript
\cite{AFV}, provisionally entitled ``Sophistication Revisited'',
the following results have been obtained. For the
set of partial recursive functions the minimal sufficient statistic
has complexity $\eqa 0$ for all data $x$. One can define equivalents of
the implicit and explicit description format in the total recursive
function setting. We obtain various upper and lower bounds
on the complexities of the minimal sufficient
statistic in all three description formats.
The complexity
of the minimal sufficient statistic for $x$, in the 
model class of total recursive functions, is called its ``sophistication.'' 
Hence, one can distinguish three different sophistications 
corresponding to the three different description formats: explicit, implicit,
and unrestricted.
It turns out that the sophistication functions are not recursive;
the Kolmogorov prefix complexity can be computed from the minimal
sufficient statistic (every description format) and {\em vice versa};
given the minimal sufficient statistic as a function of $x$ one can
solve the so-called ``halting problem'' \cite{LiVi97};
and the sophistication functions are upper semicomputable.
By the same proofs, such computability properties also hold for 
the minimal sufficient statistics in the model classes of finite sets
and computable probability mass functions.

\subsection{Application}
Because the Kolmogorov complexity is not computable, 
an algorithmic sufficient statistic cannot be computed either.
Nonetheless, the analysis gives limits to what is achievable in
practice---like in the cases of coding theorems and channel capacities
under different noise models in Shannon information theory.
The theoretical notion of algorithmic sufficient statistic forms
the inspiration to develop applied models that can be viewed
as computable approximations. 
Minimum description length (MDL),\cite{BRY}, is a good example;
its relation with the  algorithmic minimal sufficient
statistic is given in~\cite{ViLi99}. 
As in the case of ordinary probabilistic statistic,
algorithmic statistic if applied
unrestrained cannot give much insight into the {\em meaning}
of the data; in practice one must use background information to
determine the appropriate model class first---establishing what
meaning the data can have---and only then apply
algorithmic statistic to obtain the best model in that class
by optimizing its parameters. See Example~\ref{ex.restricted}.
Nonetheless, 
in applications one can 
sometimes still unrestrictedly use compression properties
for model selection, for example by a judicious
choice of model
parameter to optimize. One example is the
precision at which we represent the other parameters: too high
precision causes accidental noise to be modeled as well, too low
precision may cause models that should be distinct to be confusing.
In general, the performance of a model for a
given data sample depends critically
on what we may call the ``degree of discretization'' or the
``granularity'' of the model:
the choice of precision of the parameters, the number of
nodes in the hidden layer of a neural network, and so on.
The granularity is often determined ad hoc.
In~\cite{GLV00}, in two quite different experimental 
settings the best model granularity
values predicted by MDL are shown to coincide 
with the best values found experimentally.

\section*{Acknowledgement}
PG is grateful to Leonid Levin for some enlightening discussions
on Kolmogorov's ``structure function''.

\newpage
\begin{biography}{P\'eter G\'acs}
obtained a Masters degree in Mathematics at E\"otv\"os
University in Budapest, 1970, and a Ph.D.~in Mathematics from the J.~W.~Goethe
University, Frankfurt am Main, 1978 (advisor: C.~P.~Schnorr). He was a research
fellow at the Institute for Mathematics of the Hungarian Academy of Sciences
1970-77; a  visiting lecturer at the J.~W.~Goethe University, 1977-78;
a  research
associate in Computer Science at Stanford University in 1978-79; an
Assistent and Associate Professor at the Computer Science Department and
Mathematics Department, University of Rochester, 1980-83; and  Associate
Professor and Professor at the Computer Science Department, Boston University,
1984--present.  He was visiting scholar at Moscow State University,
University of G\"ottingen, IBM Watson Research Center, IBM Almaden Research
Center, Bell Communications Research, Rutgers University, Centrum voor
Wiskunde en Informatica (CWI).  Most of his 
papers are in the following 
areas: Shannon-type
information theory, algorithmic information theory, reliable cellular
automata.
\end{biography}

\begin{biography}{John T. Tromp}
received his Ph.D. from the University 
of Amsterdam (1993)  and he holds positions at the national 
CWI research institute in Amsterdam and is Software Developer
at Bioinformatics Solutions Inc., Waterloo, Ontario, Canada.
He has worked on cellular automata, 
computational complexity, distributed and parallel computing, 
machine learning and prediction, physics of computation, models
of computation, Kolmogorov complexity, computational biology,
and computer games.
\end{biography}

\begin{biography}{Paul M.B. Vit\'anyi}
received his Ph.D. from the Free University
of Amsterdam (1978).  He holds positions at the national
CWI research institute in Amsterdam, and he is Professor of Computer Science
at the University of Amsterdam.  He serves on the editorial boards
of Distributed Computing, Information Processing Letters,
Theory of Computing Systems, Parallel Processing Letters,
Journal of Computer and Systems Sciences (guest editor),
and elsewhere. He has worked on cellular automata,
computational complexity, distributed and parallel computing,
machine learning and prediction, physics of computation, reversible computation,
quantum computation, and algorithmic information theory
(Kolmogorov complexity). Together with Ming Li
they pioneered applications of Kolmogorov complexity
and co-authored ``An Introduction to Kolmogorov Complexity
and its Applications,'' Springer-Verlag, New York, 1993 (2nd Edition 1997),
parts of which have been translated into Chinese,
Russian and Japanese.

\end{biography}

\end{document}